\documentclass[preprint]{elsarticle}
\biboptions{sort&compress}  

\RequirePackage{amsthm,amsmath,amsfonts,amssymb}

\usepackage{bbm}
\usepackage{latexsym}
\usepackage{graphicx}
\usepackage{color,wrapfig}
\usepackage{multirow}
\usepackage{bm,url,comment}
\usepackage{mathtools,mathrsfs}
\usepackage{cool}
\usepackage{bm}
\usepackage{array}
\usepackage{colortbl}
\usepackage{adjustbox}
\usepackage{subcaption}
\usepackage{placeins}
\usepackage{setspace}
\usepackage{lmodern}
\usepackage{algorithm}
\usepackage{algpseudocode}
\usepackage{float}
\algnewcommand{\Input}{\item[\textbf{Input:}]}
\algnewcommand{\Output}{\item[\textbf{Output:}]}

\DeclareMathOperator{\OO}{O}

\theoremstyle{plain}

\theoremstyle{remark}

\newcommand{\norm}[1]{\left\lVert#1\right\rVert}

\theoremstyle{plain} 
\newtheorem{theorem}{Theorem}[section]
\newtheorem*{theorem*}{Theorem}
\newtheorem{lemma}[theorem]{Lemma}

\newtheorem*{lemma*}{Lemma}

\newtheorem*{corollary*}{Corollary}

\newtheorem*{proposition*}{Proposition}

\newtheorem{definition}[theorem]{Definition}
\newtheorem*{definition*}{Definition}
\theoremstyle{remark}

\newtheorem*{example*}{Example}
\newtheorem{remark}[theorem]{Remark}

\newtheorem*{remark*}{Remark}
\newtheorem*{remarks*}{Remarks}

\DeclareMathOperator*{\argmin}{arg\,min}
\newcommand{\lightrule}{\arrayrulecolor[gray]{0.8}\hline\arrayrulecolor{black}}  

\begin{document}

\begin{frontmatter}

\title{
 Extracting Dual Analytic Geometries of Linear Transformations to Achieve Efficient Computation}

\author[A]{Pei-Chun Su\corref{mycorrespondingauthor}}
\cortext[mycorrespondingauthor]{Pei-Chun Su}
\ead{pei-chun.su@yale.edu}
\author[A]{Ronald R. Coifman}
\address[A]{Department of Applied and Computational Mathematics, Yale University, CT, USA}

\begin{abstract}
We generalize and extend; the analytic and geometric methodology of Calder\'on--Zygmund for operator analysis and, the fast multipole computational algorithms of Rokhlin-Greengard, in which, both domain and range of an operator, are jointly organized. Our methods enable the analysis of unstructured kernels by uncovering the joint hidden geometries in their domain and range structures. A main tool is the \texttt{Questionnaire} algorithm, an iterative dynamic procedure that constructs adaptive bi-hierarchical partition trees to reveal  latent multiscale organizations and expose local low-rank structures within the kernel. When the geometry is unknown, the iteration starts from any initial affinity, even a random one, and, combined with a stochastic flip perturbation and a spectral space-filling curve that fixes the left--right order of every tree, reorganizes the matrix into an image of the operator's "phase space": the joint space of its sources and sensors, in duality, as positions and frequencies are for the Fourier transform. The ordering gives this phase space an analytic structure on which any classical transform can be carried out, and the compression it yields is the direct evidence of this structure. 
Guided by these geometries, we employ two different approaches to enable efficient matrix factorization and fast matrix vector multiplication: (1) The \texttt{Butterfly} algorithm, which exploits the learned bi-hierarchical localized low-rank structure; and (2)  Adaptive \texttt{eGHWT}, best tilings in both space and sequency using all levels of  generalized Haar--Walsh wavelet packets. 
   We coin our algorithms as \texttt{Questionnaire Factorization and Fast Transform (QFFT)}. 
 Unlike classical approaches that rely on structural prior knowledge of the underlying geometry, \texttt{QFFT} is fully data-driven and applicable to matrices arising from irregular or unknown distributions.
 We demonstrate the effectiveness of our approach on matrices associated with heterogeneous oscillatory operators, or for orthogonal eigenvector expansions in the presence of high multiplicity and random sampling. The resulting compressed representations reduce storage complexity from $\OO(N^2)$ to $\OO(N \log N)$, enabling fast computation and scalable implementation. 
\end{abstract} 
\begin{keyword}
 \texttt{QFFT};
 {Iterative Geometry Learning};
 {Kernel Factorization};
 {Fast Transform}. 
\end{keyword}

\end{frontmatter}

\section{Introduction}
A fundamental task in scientific computing is the evaluation of structured matrix–vector products of the form
\begin{equation}
u(y) = \sum_{x \in X} \mathcal{K}(y, x) f(x),
\end{equation}
where \( \mathcal{K}(y, x) \) is a kernel function defined in a discrete or continuous domain, and \( X \) denotes the set of input points. Here \( f\colon X \to \mathbb{R} \) is the input function (equivalently an input vector indexed by \( X \)) to which the operator \( \mathcal{K} \) is applied, and \( u \) is the resulting output over the target points \( Y \), so the columns of \( \mathcal{K} \) are indexed by the operator's \emph{domain} \( X \) and its rows by the \emph{range} \( Y \); we write the product as \( \mathcal{K} f \). When the kernel carries known analytic or geometric structure with respect to coordinates on its domain and range, such products can often be evaluated in $\OO(N\log N)$ rather than $\OO(N^2)$ operations by fast algorithms such as the Fast Fourier Transform, the fast multipole method (FMM)~\cite{greengard1987fast}, and the butterfly factorization~\cite{michielssen1996multilevel,candes2009fast,o2010algorithm}, all rooted in the analytic and geometric operator theory of Calder\'on--Zygmund~\cite{calderon1952existence}. Each of these methods, however, presupposes that the multiscale geometry of the rows and columns is already known and reflected in the ordering of the matrix. 

We extend the reach of these analytic and geometric decompositions to the setting in which the domain and range geometries are themselves unknown, irregular, or scrambled: we present data-driven algorithms that discover both geometries jointly and directly from the matrix entries, allowing the data to reveal natural
bi-structural analogs.
\subsection{Problem Statement}

We consider a single problem in two guises. Given a linear operator presented only as an $N_Y \times N_X$ matrix $\mathcal{K}$, with its rows and columns in unknown or arbitrary order and no accompanying coordinates, our goal is to develop a fast matrix--vector product algorithm $\mathcal{K} f$, rather than the original $\OO(N^2)$ of a dense multiplication, $\mathcal{K}$ exhibits no organized structure and no fast algorithm applies to it directly. This does not imply that no such structure exists: for many operators, such as the  Fourier Transform or a discrete sine transform, a multi scale coordinates ordering admitting an $\OO(N\log N)$ algorithm \emph{is} present, but is unknown to us in the absence of coordinates on the domain and range. The task is therefore to discover, directly from the matrix entries and without prior knowledge of the underlying domain-range  geometry, a multiscale organization of the rows and columns under which such a fast product exists.

We model the observed matrix as a row/column permutation of an analytically structured one. We consider multiscale partition trees as basic building blocks.

\begin{definition}[Multiscale Partition Tree Structure \( T_X \)]\label{def:partition_tree}
Let \( X = \{x_i\}_{i=1}^N \) be a set of sample points. A multiscale partition tree \( T_X \) is a hierarchical decomposition of the index set \( \{1, \dots, N\} \) into disjoint subsets, referred to as \emph{folders}, denoted $\left\{ V_k^\ell \right\}_{0 \le \ell \le L,\, 0 \le k < 2^\ell }$, where \( \ell = 0, \dots, L \) indexes the level, \( k \) indexes the folders at each level \( \ell \), and \( \mathcal{C}_k^\ell \) denotes the index set of the child folders of \( V_k^\ell \). The root of the tree is \( V_0^0 = \{1, \dots, N\} \), and each folder \( V_k^\ell \) at level \( \ell \) is the disjoint union of its children at level \( \ell+1 \):
\begin{equation}
V_k^\ell = \bigcup_{m \in \mathcal{C}_k^\ell} V_m^{\ell+1}.
\end{equation}
\end{definition}

In this work we limit the discussion to \emph{balanced} binary trees. Such a tree, together with a left--right ordering of the two children of every folder, defines a permutation of the points; Section~\ref{sec:sfc} makes this precise.

Writing $\mathfrak{S}_n$ for the symmetric group of permutations of $\{1,\dots,n\}$, we assume there exist a column permutation $\pi_X \in \mathfrak{S}_{N_X}$ and a row permutation $\pi_Y \in \mathfrak{S}_{N_Y}$ for which the reordered matrix
\begin{equation}\label{eq:problem}
  \mathcal{K}^{\pi}_{ji} \;:=\; \mathcal{K}\!\left(\pi_Y(j),\,\pi_X(i)\right),
  \qquad 1 \le j \le N_Y,\;\; 1 \le i \le N_X,
\end{equation}
is compressible with respect to the tree structure and left--right ordering that the two permutations define: read in the orders $\pi_X$ and $\pi_Y$, the folders of balanced binary trees on the columns and rows are dyadic intervals, and the blocks of $\mathcal{K}^{\pi}$ over pairs of such folders are of low rank across dyadic scales, so that $\mathcal{K}^{\pi}$ is described by $\OO(N\log N)$ parameters and admits fast algorithms that exploit this bi-hierarchical structure. Given only the entries of $\mathcal{K}$, our goal is to find such a pair $(\pi_X,\pi_Y)$; throughout we take $N_X,N_Y = \OO(N)$. We build on two fast algorithms of this kind.

\texttt{Butterfly} factorization skeletonizes well-separated sub-blocks so that the reordered kernel telescopes into a product $\mathcal{K}^{\pi}\approx B^{(L)}P^{(L-1)}\cdots P^{(0)}$ of a coarse skeleton block $B^{(L)}$ and a chain of interpolation matrices $P^{(\ell)}$ (\ref{app_Butterfly}). 

The \texttt{eGHWT} instead selects a sparse wavelet-packet best basis,  expanding the kernel over two-dimensional Walsh atoms $w_{(q,p)}(y,x)=w_q(y)\,w_p(x)$ and choosing the tiling $\mathcal{B}^{*}$ of least coefficient cost (\ref{app_Walsh}).

This formulation covers the two settings we study:
\begin{enumerate}
\item Kernels $\mathcal{K}(y,x)$ arising from oscillatory or potential-type interactions between sources $x$ and sensors $y$, whose bi-hierarchical dual structure is hidden by an unknown ordering.
\item Empirical eigenbases: eigenfunctions of the Laplacian on a manifold sampled at random points, either known in closed form (the spherical harmonics, Section~\ref{sec:sh}) or computed as eigenvectors of a graph Laplacian (Section~\ref{sec:manifold}). For such samples no classical fast transform is available, in particular under high spectral multiplicity. The points and the eigenvectors are in duality, like positions and frequencies, so the reorganized matrix of eigenvectors is a picture of the phase space of the manifold (Section~\ref{sec:sfc}).
\end{enumerate}
 
\subsection{Unsupervised Ground Metric Learning}\label{sec:related}

The problem posed above is related to  \emph{unsupervised ground metric learning}: from the data matrix alone, with no labels or reference coordinates, one must recover a coupled row--column geometry, a distance on each axis.

We use these distances to reorder  the rows and columns of the matrix into a structured  compressible form revealing their dual hidden structures. 
Our goal is to  learn it directly from the data, without supervision, by exploiting a coupling between the two modes of a data matrix: A good geometry on the samples (rows) informs the geometry on the features (columns), and vice versa. 
Classically this metric is supplied through labels or side information, where it is termed a \emph{ground metric}~\cite{cuturi2014ground}; the methods reviewed below, ours among them, instead learn it, together with the coupled geometries, it encodes, from the matrix alone.

This coupled, alternating view was introduced by the \texttt{Questionnaire} algorithm~\cite{coifman2011harmonic,gavish2012sampling,ankenman2014geometry}, the geometry-learning engine of this work. The name reflects its motivating example, a survey questionnaire whose rows are questions and columns are respondents, for which organizing the data means simultaneously grouping the questions that elicit similar answers and the respondents who answer alike, the dual row--column geometry of the algorithm is inferred from the response patterns alone. It organizes the rows and columns of a matrix into mutually informed multiscale partition trees: a tree on one axis induces a dual affinity on the other, from which a new tree is built, and the two geometries are refined in alternation (until convergence) .
A recent variant in the same lineage~\cite{lin2026joint} modifies only the tree weights of the dual geometry, reading them from a hyperbolic embedding~\cite{nickel2017poincare} of the tree.
We review the \texttt{Questionnaire} algorithm in detail in Section\ref{sec:background}. In the examples reported in prior work, the iteration is seeded with an informative initial geometry, typically a cosine-similarity or Euclidean-distance affinity computed from the data. In Section~\ref{sec:dstiv} we instead demonstrate it on the discrete sine transform of type IV (DST-IV), whose rows and columns are mutually orthogonal, so no such informative geometry is available; there we start instead from a Gaussian random affinity $W_X^0$ on the columns and let the coupled iteration discover structures from the matrix alone.  

\subsection{Demonstration: The Scrambled DST-IV}\label{sec:dstiv}

   We now illustrate the approach on a concrete, deliberately adversarial instance of this problem: the DST-IV, whose kernel is $\sin\!\big(\pi(y-\tfrac12)\,x\big)$ at integer frequencies $y\in\{1,\dots,N\}$ and spatial sample points $x\in\big\{\tfrac{n-1/2}{N}:n=1,\dots,N\big\}$, with $N=256$. We hide the natural ordering: let $\pi_Y,\pi_X\in\mathfrak{S}_N$ be two independent, uniformly random permutations, so that observed row $j$ carries the frequency $y_j=\pi_Y(j)$ and observed column $i$ the spatial coordinate $x_i=\tfrac{\pi_X(i)-1/2}{N}$. The observed matrix is then

\begin{equation}
\mathcal{K}_{ji} = \mathcal{K}(y_j, x_i) = \sin\!\big( \pi (y_j-\tfrac12)\, x_i \big).
\end{equation}
  
   Only the scrambled entries $\mathcal{K}_{ji}$ are given (neither the permutations $\pi_X,\pi_Y$ nor the fact that a DST-IV underlies them), and discovering the hidden frequency and space orders is exactly the problem posed above.

\begin{figure}[t]
\centering
\includegraphics[width=\linewidth]{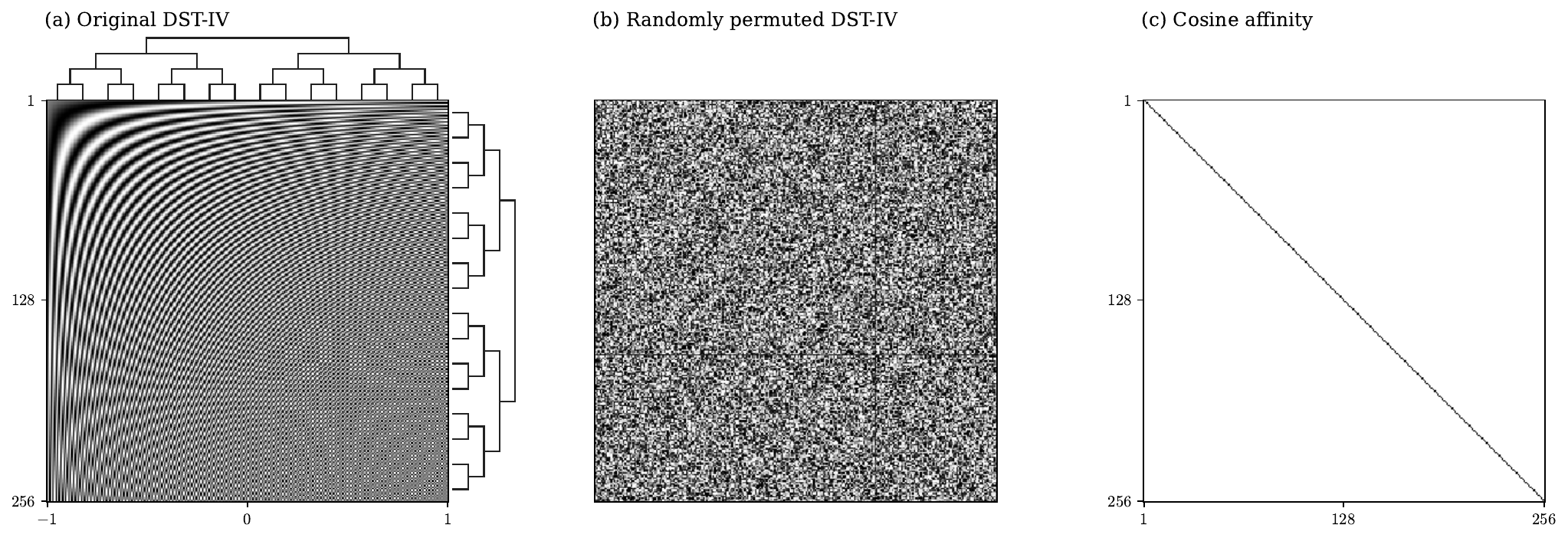}
\caption{The scrambled DST-IV problem ($N=256$): \textbf{(a)} natural order with the hidden row and column trees, \textbf{(b)} the observed random permutation, \textbf{(c)} the cosine or metric affinity.}
\label{fig:problem_dstiv}
\end{figure}

In its scrambled form the matrix exhibits no exploitable structure: applying it to a vector costs $\OO(N^2)$, and no fast algorithm is available. Figure~\ref{fig:problem_dstiv} illustrates the setup. Panel~(a) shows the DST-IV in its natural order, smoothly organized and carrying the hidden row and column partition trees; panel~(b) shows the randomly permuted matrix we actually observe, with no visible structure. Panel~(c) shows the cosine affinity in the observed order: because the DST-IV has orthonormal rows \emph{and} columns, this affinity is the identity on \emph{both} sides, so neither the range nor the domain exposes any geometry to start from.

\begin{figure}[p]
    \centering
    \begin{subfigure}[b]{\textwidth}
        \centering
        \makebox[\textwidth][c]{\includegraphics[width=\linewidth,trim=3cm 1cm 3cm 1cm, clip]{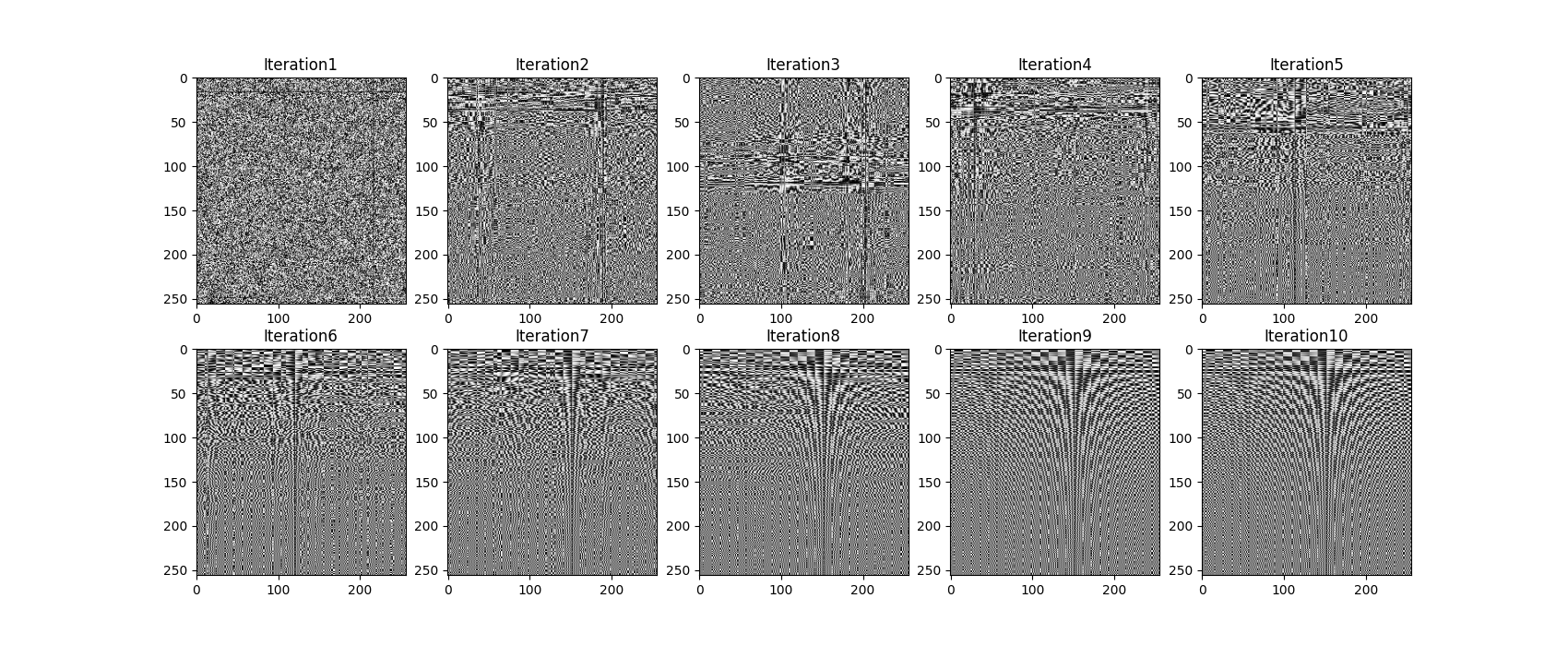}}
        \caption{ Reorganization.}
        \label{fig:sin_reorg}
    \end{subfigure}

    \vspace{0.4cm}

    \begin{subfigure}[b]{\textwidth}
        \centering
        \makebox[\textwidth][c]{%
            \includegraphics[width=\linewidth, trim=2cm 2cm 3cm 1cm, clip]{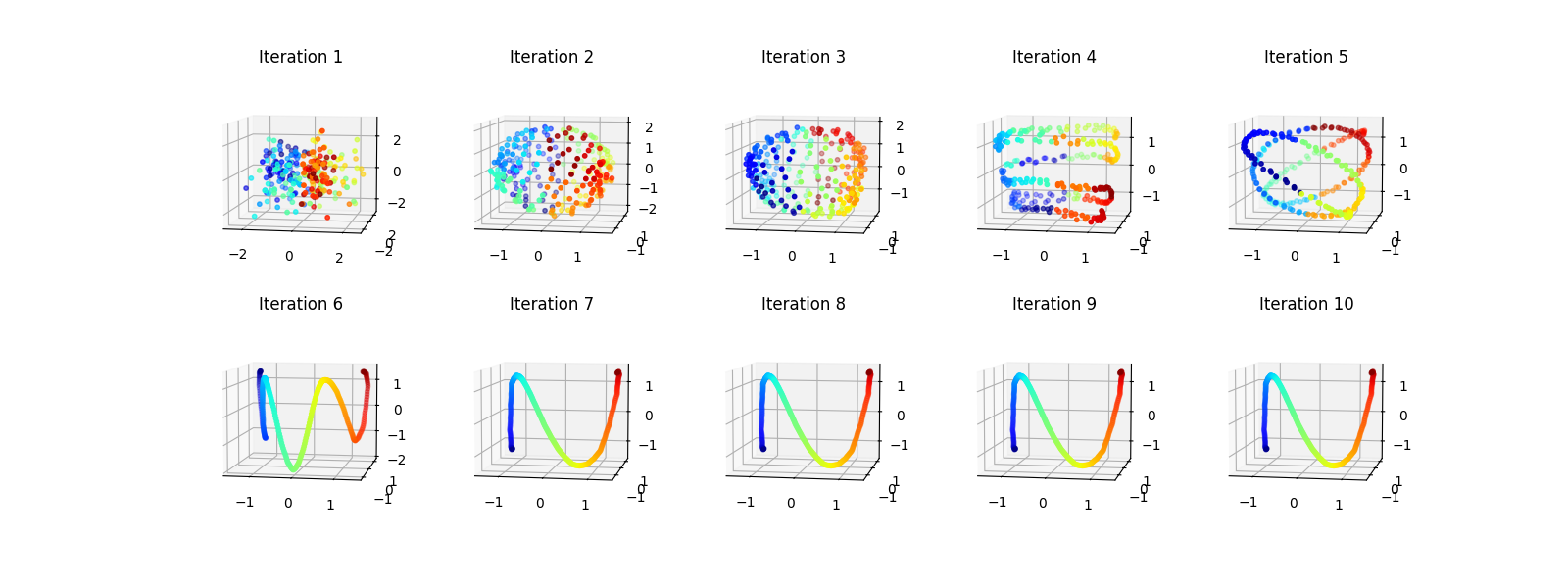}%
        }
        \caption{Embedding of the space $X$.}
        \label{fig:sin_col}
    \end{subfigure}

    \vspace{0.5cm}

    \begin{subfigure}[b]{\textwidth}
        \centering
        \makebox[\textwidth][c]{%
            \includegraphics[width=\linewidth, trim=2cm 2cm 3cm 1cm, clip]{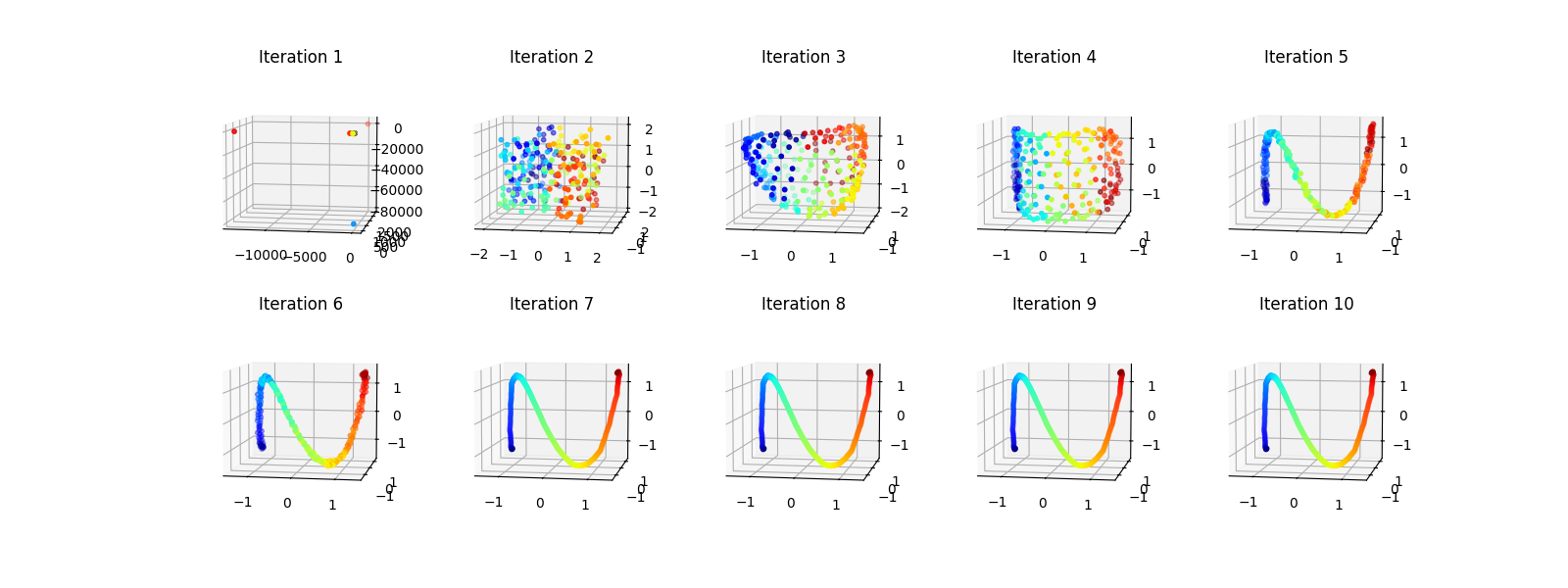}%
        }
        \caption{ Embedding of frequency space $Y$.}
        \label{fig:sin_row}
    \end{subfigure}

    \caption{Analysis of the "scrambled" DST-IV kernel over 10 iterations of Algorithm~\ref{alg:questionnaire_basic}, with rows and columns reordered by the respective space-filling curve of Algorithm~\ref{alg:affinity_ordering}: showing (a) reorganization, (b) space embedding, and (c) frequency embedding.}
    \label{fig:sin_combined_vertical}
\end{figure}

Because this affinity offers no starting geometry, we seed the column iteration with a Gaussian random affinity: we draw $N$ points $z_i$ i.i.d.\ from a standard Gaussian in $\mathbb{R}^{10}$ and set $W_X^0(i,i') = \exp\!\big(-\lVert z_i - z_{i'}\rVert^2/\sigma^2\big)$, with $\sigma$ the median pairwise distance. We then run 10 iterations of the \texttt{Questionnaire} algorithm (Algorithm~\ref{alg:questionnaire_basic}). The dual affinity is the multiscale correlation affinity (Definition~\ref{def_Tree_Aff}, Equation~\eqref{eq:corr}), from which the row and column trees and their leaf ordering are constructed at each iteration by Algorithm~\ref{alg:affinity_ordering}, gradually rearranging the matrix (Figure~\ref{fig:sin_combined_vertical}).
Although the rows and columns of the matrix are orthogonal and randomly permuted, our method successfully reveals the intrinsic local joint structure of the phase space of frequencies (rows) and positions (columns), reorganizing the matrix into a highly structured and symmetric form as the iterations proceed.
The diffusion embeddings of the newly sorted rows and columns show that the method gradually uncovers the underlying one-dimensional manifold structure of the spatial domain $X$ as well as the frequency domain $Y$ (Figures~\ref{fig:sin_col} and~\ref{fig:sin_row}). The smooth color gradients along the embeddings reflect the specific reorderings, illustrating how hidden geometric relationships in space and frequency (phase space) are progressively and successfully discovered.

\begin{figure}[t]
\centering
\includegraphics[width=1\linewidth]{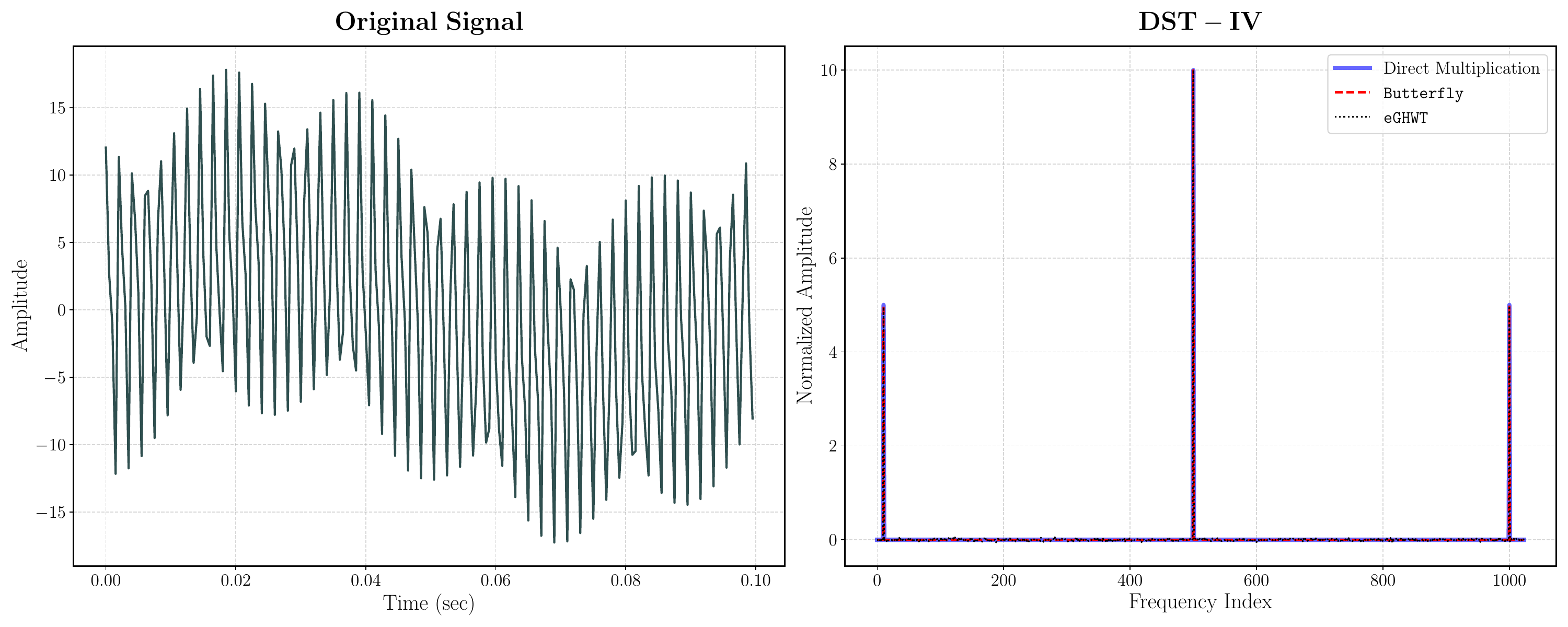
}
\caption{Approximation of the DST-IV transform ($N = 256$) using \texttt{Butterfly} and \texttt{eGHWT}. \textbf{Left}: the input signal
\(f(x) = 5\sin(\pi (y_1+0.5)x) +
          10\sin( \pi (y_2+0.5)x) +
          5\sin( \pi (y_3+0.5)x))\), where $y_1 = 10$, $y_2 = 100$, and $y_3 = 200$.
\textbf{Right}: the transform by direct multiplication, by \texttt{Butterfly} (target accuracy $\varepsilon = 10^{-6}$, $\max_k|V_k^\ell| = 32$, measured relative $\ell_2$ error $1.6\times10^{-6}$), and by \texttt{eGHWT} (target accuracy $\varepsilon = 5\times 10^{-2}$, measured relative $\ell_2$ error $5.0\times10^{-2}$).}
\label{fig:sin_transform}
\end{figure}

Once the hidden row and column orders are discovered, the matrix is revealed to be a (permuted) discrete sine transform, which admits an $\OO(N\log N)$ fast algorithm. In effect, our procedure \emph{rediscovers a fast transform directly from the scrambled data}, with no prior knowledge that a DST-IV structure was present. To confirm that the rediscovered ordering delivers genuine fast transforms with high precision, we use it to apply the operator. With the DST-IV kernel above at $N = 256$ and the test signal
\begin{equation}
f(x_i) = 5\sin(\pi (y_1+0.5)x_i) +
          10\sin( \pi (y_2+0.5)x_i) +
          5\sin( \pi (y_3+0.5)x_i),
\end{equation}
where $y_1 = 10$, $y_2 = 100$, and $y_3 = 200$, we approximate the transform $\int_{0}^{1} \mathcal{K}(y,x)\, f(x)\, dx$ by applying the reorganized kernel through the fast \texttt{Butterfly} or the two-dimensional \texttt{eGHWT} factorizations (the fast-application derivations are given in Section~\ref{section_kerenl_app} and \ref{app_Walsh}). Figure~\ref{fig:sin_transform} shows the result: on the left the input signal, and on the right its DST-IV transform, which concentrates the signal's energy into three sharp spikes at exactly the input frequencies $y_1 = 10$, $y_2 = 100$, and $y_3 = 200$, recovering their amplitudes $5$, $10$, and $5$. The transforms computed by direct dense multiplication, by \texttt{Butterfly}, and by \texttt{eGHWT} coincide on these peaks, confirming that the representation learned from the scrambled matrix accurately reproduces the action of the transform. Each factorization is built to a prescribed accuracy $\varepsilon$, and the measured relative $\ell_2$ error of the resulting fast matrix-vector product matches it closely: the \texttt{Butterfly} at $\varepsilon = 10^{-6}$ (maximum leaf-block size $\max_k|V_k^\ell| = 32$; see~\ref{app_Butterfly}) attains $1.6\times10^{-6}$, and the \texttt{eGHWT} at $\varepsilon = 5\times10^{-2}$ attains $5.0\times10^{-2}$.

Table~\ref{tab:timing} compares the single-thread computation time of each algorithm on this example. Preparation dominates the cost. Each \texttt{Questionnaire} iteration takes about $141$ ms and scales as $\OO(N^3)$: it forms the two $N\times N$ dual affinities as pairwise tree--earth-mover distances and rebuilds the trees by spectral bisection. Factorizing the reorganized kernel is a further one-time cost that reads the full $N\times N$ matrix: $520$ ms for the \texttt{eGHWT}, whose best-basis selection scales as $\OO(N^2(\log N)^2)$, and $9.5$ ms for the \texttt{Butterfly}, whose construction scales as $\OO(N^2\log N)$. Once built, the operator applies to each vector in only $\OO(N\log N)$: $0.07$ ms per matrix--vector product for the \texttt{Butterfly}, and $14$ ms for the \texttt{eGHWT}, which shares the $\OO(N\log N)$ application scaling but remains slower in absolute terms in our single-threaded implementation. The one-time preparation is thus computationally intensive, but it is precisely what \emph{discovers} a fast $\OO(N\log N)$ algorithm for an operator that admitted none in its scrambled ordering, so it is paid only once, after which every application of the transform is cheap.

\begin{table}[t]
\centering
\begin{tabular}{lr}
\hline
Stage (single thread, $N=256$) & Time \\
\hline
\texttt{Questionnaire}, per iteration & $141$ ms \\
\texttt{eGHWT} preparation & $520$ ms \\
\texttt{Butterfly} preparation & $9.5$ ms \\
\texttt{eGHWT} multiplication (per vector) & $14$ ms \\
\texttt{Butterfly} multiplication (per vector) & $0.07$ ms \\
\hline
\end{tabular}
\caption{Average single-thread computation times on the scrambled DST-IV ($N=256$): one \texttt{Questionnaire} iteration; the \texttt{eGHWT} and \texttt{Butterfly} preparation (factorization of the Reorganization-2 iterate, with the same accuracy parameters as Figure~\ref{fig:sin_transform}); and the operator application, averaged over $100$ random vectors.}
\label{tab:timing}
\end{table}

\begin{figure}[t]
\centering
\includegraphics[width=\linewidth]{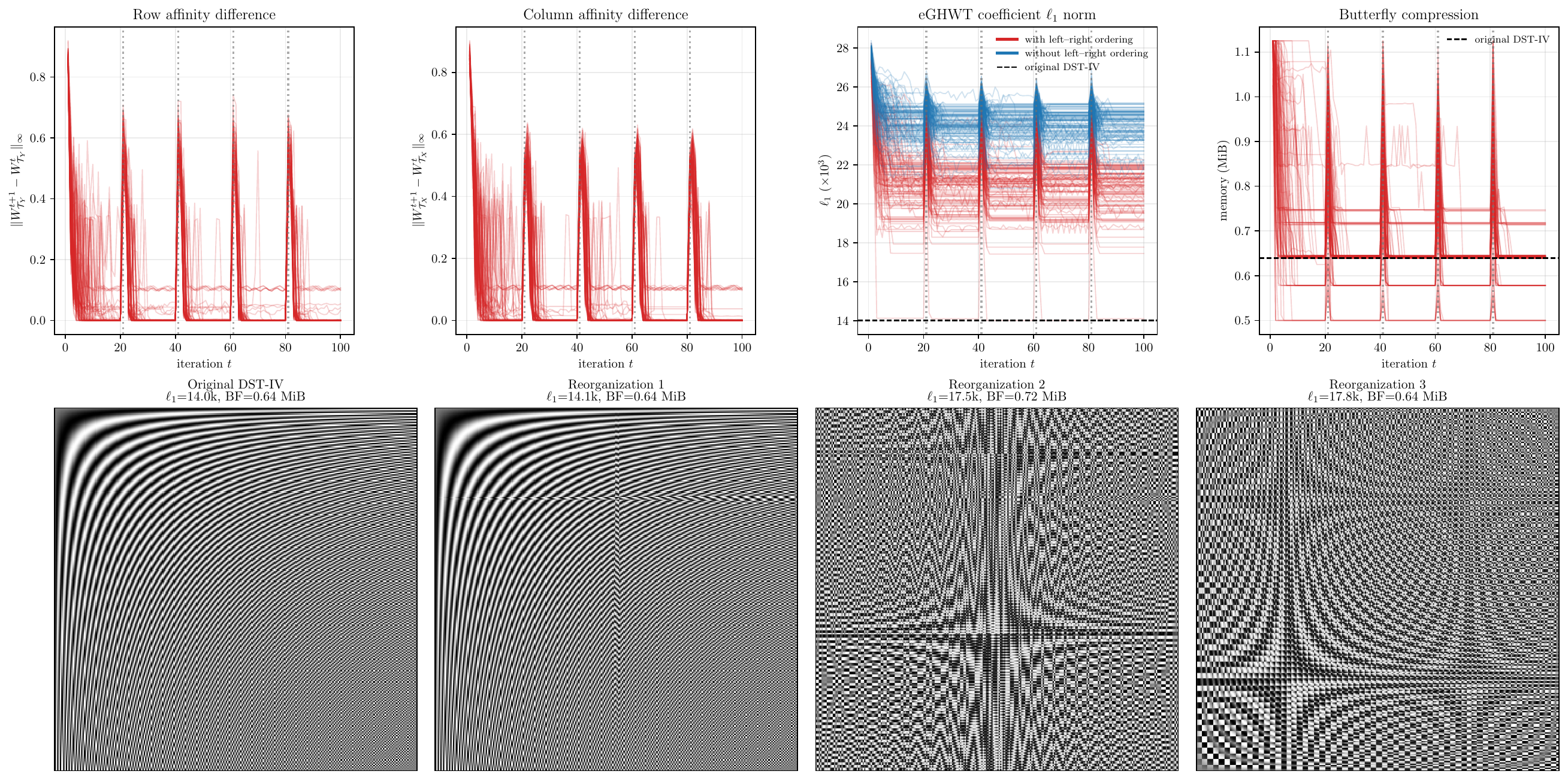}
\caption{The \texttt{Questionnaire} with the flip perturbation on the scrambled DST-IV ($N=256$), from $100$ random initializations $W_X^0$, with a perturbation applied every $20$ iterations (dotted vertical lines). \textbf{Top}, left to right: the row and column affinity differences $\|W_Y^{t+1}-W_Y^{t}\|_\infty$ and $\|W_X^{t+1}-W_X^{t}\|_\infty$, the $\ell_1$ norm of the \texttt{eGHWT} best-basis coefficients (at $\varepsilon=5\times10^{-2}$) with (red) and without (blue) the left--right ordering of the space-filling curve, and the \texttt{Butterfly} memory (at $\varepsilon=10^{-6}$, $\max_k|V_k^\ell|=32$), versus iteration $t$; dashed lines mark the natural DST-IV. \textbf{Bottom}: the original DST-IV and three reorganizations at the final iteration, each annotated with the $\ell_1$ norm of its \texttt{eGHWT} coefficients and \texttt{Butterfly} memory.}
\label{fig:q100_fixedpoints}
\end{figure}

The reorganization in Figure~\ref{fig:sin_combined_vertical} is one \emph{good} fixed point, but not the only limit the raw iteration can reach. Since the \texttt{Questionnaire} map is discontinuous and the seed $W_X^0$ is random, different initializations behave very differently. Figure~\ref{fig:q100_fixedpoints} runs the \texttt{Questionnaire} with the flip perturbation (Algorithm~\ref{alg:questionnaire}) from $100$ random seeds. Over the first $20$ iterations, which coincide with the raw \texttt{Questionnaire} (Algorithm~\ref{alg:questionnaire_basic}) before any perturbation is applied, the row and column affinity differences mostly decay, but the limit the iteration reaches (a fixed point for some seeds, a short periodic cycle for others) varies widely in quality: measured by either the $\ell_1$ norm of the \texttt{eGHWT} best-basis coefficients (at $\varepsilon = 5\times10^{-2}$, \ref{app_Walsh}) or the \texttt{Butterfly} memory (at $\varepsilon = 10^{-6}$, $\max_k|V_k^\ell| = 32$), a few land close to the natural DST-IV (dashed) while most sit well above it. The \texttt{eGHWT} panel further shows this $\ell_1$ norm computed with (red) and without (blue) the left--right ordering of the tree. The \texttt{eGHWT} best basis is built on a hierarchical partition tree, but its cost depends on the left--right ordering of the tree's sibling nodes: this ordering fixes how the leaves are laid out on the line, and hence the smoothness of the reorganized matrix and the sparsity of the best basis. It is left unspecified by the original construction~\cite{saito2022eghwt} and is not unique in coordinate-free settings; our space-filling curve (Algorithm~\ref{alg:affinity_ordering}) supplies a principled choice. Imposing it (red) consistently lowers the $\ell_1$ norm of the same organization relative to leaving it unordered (blue), so the left--right ordering is itself a source of compressibility for the \texttt{eGHWT}.

Because the reachable limits vary so widely in quality, monitoring the change of the affinity certifies neither convergence nor the quality of the organization; only the transform cost does. This motivates the two ingredients we add to the raw iteration: tracking the compression cost directly, to select the most compressible iterate, and a stochastic \emph{flip} perturbation that dislodges a weak fixed point or periodic cycle and lets the iteration drop to one whose reorganization is smoother, and hence more compressible, with respect to the learned dual tree geometry. In Figure~\ref{fig:q100_fixedpoints} this perturbation, detailed in Section~\ref{subsec:perturbation}, is applied every $20$ iterations, each time holding a random half of the points on each axis fixed and permuting the other half, producing the sawtooth beyond iteration~$20$. Each perturbation dislodges the current limit and the iteration re-descends toward a stronger one; by the final step the three reorganized matrices at the bottom all recover coherent, locally low-rank block structure rather than the random scrambling of the input.

\subsection{Our Contributions}
 We enable ideas for processing and structuring data bases \texttt{Questionnaire} algorithm~\cite{coifman2011harmonic} (Algorithm~\ref{alg:questionnaire_basic}) to organize the analysis of complex seemingly unstructured operators 
We discover the hidden multiscale structure without predefined coordinates, to restructure a linear transformation,and build adaptive, multiscale \emph{ordered} binary trees jointly on the rows and columns, iteratively, reducing variation across scales. Each tree, together with its left--right ordering, is constructed by Algorithm~\ref{alg:affinity_ordering} (Section~\ref{sec:sfc}). Through the resulting tree-based metrics, the algorithm finds organized geometries of the kernel, the fixed points, or cycles, of the tree-building iteration, and yields a multiscale basis tailored to the kernel values, that is, to the interactions between inputs and outputs rather than to the physical coordinates of the domain and range.

To this iteration we add two ingredients, which together let it start from any initial affinity, even a random one, and still recover a compressed "analytic" structure of the operator's phase space. First, from the mutually learned trees we construct two data-driven space-filling curves (Algorithm~\ref{alg:affinity_ordering}), one aligned with the intrinsic geometry of the rows and one with that of the columns. Algorithm~\ref{alg:affinity_ordering} fixes the left--right order of the sibling subtrees as part of building each tree, an order that bottom-up constructions leave arbitrary, and the resulting curves align the rows and columns out along the learned geometry. The reorganized matrix then becomes an "image" of the operator's "phase space", the joint space of its sources (the domain $X$) and receivers (the range $Y$), in which the two are in duality as positions and frequencies are for the Fourier transform; Section~\ref{sec:sfc} makes this phase-space picture precise. This dual ordering  gives the phase space an analytic structure. The reorganized matrix becomes a function of the two curve coordinates (Section~\ref{sec:sfc}), enabling tools of one- and two-dimensional analysis, from Fourier and local trigonometric bases to wavelets, wavelet packets and Walsh tilings,to apply  directly, any transform defined on ordered data can be carried out on the learned organization. The \texttt{Best basis Tensor Walsh packet} transform,eGHWT is an the example we  study: it requires a left--right ordering of the trees. The compression gain is the direct evidence of this structure: for the same learned organization. Imposing our ordering measurably lowers the \texttt{eGHWT} cost (Figure~\ref{fig:q100_fixedpoints}). Second, because the iteration can stall at a weak fixed point or a periodic cycle, in particular from an uninformative random start, a flip perturbation (Algorithm~\ref{alg:questionnaire}, Section~\ref{subsec:perturbation}) dislodges that limit and continues the iteration toward a more compressible organization.

This organization is precisely what the two compression methods ,that we use, require: \texttt{Butterfly} factorization exploits localized low-rank interaction blocks, and the \texttt{eGHWT} a sparse wavelet-packet best basis, both of which emerge only once the rows and columns are organized into geometrically coherent sub-blocks, and which in the classical setting require the known coordinates unavailable here. From geometry discovery to structured factorization, this pipeline defines our \emph{Questionnaire Factorization and Fast Transform} (\texttt{QFFT}) algorithm.

In particular, the reorganized matrix can be compressed as if it were an image. Equivalently, the method implicitly ``dilates'' the operator on $\mathbb{R}^N$ into a structured sparse matrix in a higher-dimensional space $\mathbb{R}^{N\log N}$, in which the non-standard matrix--vector product costs only $\OO(N\log N)$. This non-standard matrix--vector multiplication is realized by expanding a vector in $\mathbb{R}^N$ in the standard Walsh basis, computing $N\log N$ numbers along the binary tree as the coordinates  of the lifted vector on which the sparse compressed matrix is applied, followed by the transposed map \cite{beylkin1991fast}.

\subsection{Organization of the Paper}

 Section 2 provides the background on unsupervised geometry learning: spectral graph theory and the coupled row and column geometry learned via the \texttt{Questionnaire} algorithm. Section 3 presents our method, comprising the stochastic flip perturbation that escapes weak fixed points and periodic cycles, the recursive tree construction and spectral space-filling curve that reorder the rows and columns, and the resulting \texttt{QFFT} fast-multiplication algorithm. Section 4 describes our evaluation methodology, including the hardware and software specifications used for the numerical simulations, and Section 5 presents comprehensive numerical experiments that quantify the efficiency and effectiveness of the \texttt{QFFT} across various structured operators. Section 6 concludes with a discussion and directions for future work. Supporting background material is collected in the appendices: the Walsh wavelet-packet and best-basis construction, interpolative decomposition, and the \texttt{Butterfly} factorization.

\section{Background of Unsupervised Learning}\label{sec:background}

Having outlined in the introduction how a learned row--column geometry feeds the \texttt{Butterfly} and \texttt{eGHWT} compression engines, we now develop that geometry-learning framework in detail. Both methods require a multiscale tree structure on the rows and columns, which for the coordinate-free operators we target is not available in advance; we therefore construct it directly from the kernel values. The framework rests on two ingredients, which we review in turn: spectral graph theory (Section~\ref{sec:graph}), which extracts a hierarchical partition from an affinity matrix, and the \texttt{Questionnaire} algorithm, which discovers the dual dependencies between row and column affinities and their coupled multiscale structure from the matrix itself.

\subsection{Spectral Graph Theory}\label{sec:graph}

Spectral graph theory provides a powerful framework for analyzing functions defined on graphs using the eigenstructure of matrices derived from graph connectivity. Given an undirected weighted graph
$G = (V, E)$  with $N$ nodes, we define the \emph{edge weight matrix} (or adjacency matrix) $W \in \mathbb{R}^{N \times N}$ such that $W_{ij} \geq 0$ encodes the affinity or similarity between nodes $i$ and $j$, with $W_{ij} = 0$ if there is no edge between them. The \emph{degree matrix} $D \in \mathbb{R}^{N \times N}$ is a diagonal matrix with entries
$D_{ii} := \sum_{j=1}^{N} W_{ij}$,
which represents the total connection strength of node $i$ to all other nodes.
From these matrices, one defines several forms of graph Laplacians, including the unnormalized Laplacian $L = D - W$, the random-walk normalized Laplacian $L_{\text{rw}} = D^{-1}L$, and the symmetric normalized Laplacian
$L_{\text{sym}} = D^{-1/2} L D^{-1/2} = I - D^{-1/2} W D^{-1/2}$.
When the vertices $V$ are viewed as a set of $N$ points sampled i.i.d. from a compact $m$-dimensional manifold $\mathcal{M} \subset \mathbb{R}^d$, the random-walk normalized Laplacian converges to the Laplace-Beltrami operator $\Delta_{\mathcal{M}}$ \cite{singer2006graph}.  The unnormalized Laplacian $L$ and symmetric Laplacian $L_{\text{sym}}$ also converge to $\Delta_{\mathcal{M}}$ up to a density-dependent scaling factor \cite{coifman2006diffusion,nadler2005diffusion}.
These graph Laplacians play a central role in applications such as spectral clustering, graph signal processing, and dimensionality reduction.
In particular, the random-walk Laplacian is closely related to the diffusion operator used in \emph{diffusion maps} \cite{coifman2006diffusion}, a nonlinear dimensionality reduction technique that exploits the long-time behavior of random walks on graphs. Its eigenvectors form a data-driven coordinate system that preserves the geometry of the underlying data manifold. Of particular interest is the second eigenvector of $L_{\text{sym}}$, known as the \emph{Fiedler vector}, which captures the most significant nontrivial mode of variation on the graph. This vector serves as a foundational tool in our framework, guiding the construction of adaptive multiscale representations through hierarchical graph partitioning.

The Fiedler vector is a natural object for this construction because it is the continuous relaxation of a balanced graph cut. It minimizes the degree-weighted Rayleigh quotient
\begin{equation}\label{eq:fiedler}
\mathbf{f} \;=\; \argmin_{\,\mathbf{f}^{\top} D\mathbf{1}=0}\;
\frac{\sum_{(i,j)\in E} W_{ij}\,(f_i - f_j)^2}{\sum_{i} D_{ii}\,f_i^{2}}\,,
\qquad\text{equivalently}\qquad L\mathbf{f}=\lambda D\mathbf{f}\,,
\end{equation}
whose numerator is a Dirichlet energy penalizing the separation of strongly connected (high-affinity) vertices, and whose degree-weighted denominator and constraint $\mathbf{f}^{\top}D\mathbf{1}=0$ rule out the trivial constant solution. Rounding the entries to $\{\pm1\}$ makes this the (NP-hard) \emph{normalized} minimum cut, which balances the two sides by connectivity (degree mass) rather than by cardinality and is therefore robust to the uneven vertex degrees of a data affinity; relaxing to real values gives the generalized eigenproblem above, which we solve in its symmetric form $L_{\mathrm{sym}}=D^{-1/2}LD^{-1/2}$ for a stable, real spectrum. Thresholding $\mathbf{f}$ separates the vertices into two groups joined by weak affinities: by Cheeger's inequality the best threshold gives a near-optimal normalized cut, and the median threshold, which we use, gives two halves of equal size. Applied recursively, median bisection by the Fiedler vector builds a balanced partition tree, but it leaves unspecified which of the two halves at each split comes first. Algorithm~\ref{alg:affinity_ordering} (Section~\ref{sec:sfc}) performs this bisection and adds the one step it lacks: at every split it orients the two halves by the seam rule~\eqref{eq:seam}, so that the tree comes with a left--right ordering of its leaves, the space-filling curve of Section~\ref{sec:sfc}.

\subsection{Learning the Coupled Geometry by Questionnaires}\label{sec:questionnaire}
In this section, we review the framework of the \texttt{Questionnaire} algorithm introduced in~\cite{coifman2011harmonic, gavish2012sampling}, which provides a multiscale approach to uncover the underlying coupled geometry of a matrix by organizing its rows and columns iteratively.
The key advantage of this approach is that it enables the discovery of intrinsic structures without prior knowledge of the geometry of either the row or column spaces. 
By identifying the most informative multiscale partitions, it projects high-dimensional rows and columns into a lower-dimensional representation that preserves the essential structural features of the original matrix.
Treating the matrix entries as interactions between dual entities, the geometry is inferred directly from the data.
We begin with definitions of multiscale affinities.
\begin{definition}[Tree Affinities]\label{def_Tree_Aff}
Let \( X = \{x_1,\ldots,x_N\} \) be a dataset, and let $\mathcal{T}
_X= \left\{ V_k^\ell \right\}_{0 \le \ell \le L,\, 0 \le k < 2^\ell }$ denote a multiscale partition tree over \( X \), where \( V_k^\ell \) is the \( k \)-th folder at level \( \ell \). For any function \( f : X \to \mathbb{R} \), define the sample vector over folder \( V_k^\ell \) as
\begin{equation}
f(V_k^\ell) := [f(x)]_{x \in V_k^\ell}.
\end{equation}
Let the weight function \( \omega : \mathcal{T}_X \to \mathbb{R} \) be given by
\begin{equation}
\omega(V_k^\ell) := 2^{-\alpha \ell } \cdot |V_k^\ell|^{\beta},
\end{equation}
where \( \alpha \in \mathbb{R} \) controls level sensitivity, and \( \beta \in \mathbb{R} \) modulates the influence of folder size.

\begin{itemize}
\item \textbf{(Tree-based EMD)}: Given two functions \( f, g : X \to \mathbb{R} \), the multiscale Earth Mover's Distance over \( \mathcal{T}_X \) is defined as
\begin{equation}\label{eq:emd}
EMD_{\mathcal{T}_X}(f, g) := \sum_{V_k^\ell \in \mathcal{T}_X} \left\| f(V_k^\ell) - g(V_k^\ell) \right\|_{1} \cdot \frac{\omega(V_k^\ell)}{|V_k^\ell|}.
\end{equation}
The corresponding tree affinity is
\begin{equation}
W_{\mathcal{T}_X}^{\mathrm{EMD}}(f, g) := \exp\left( -\frac{EMD_{\mathcal{T}_X}(f, g)}{\epsilon} \right),
\end{equation}
where \( \epsilon > 0 \) is a scaling parameter.

\item \textbf{(Tree-based Correlation)}: The multiscale correlation affinity over \( \mathcal{T}_X \) is defined as
\begin{equation}\label{eq:corr}
W_{\mathcal{T}_X}^{\mathrm{corr}}(f, g) := \sum_{V_k^\ell \in \mathcal{T}_X} \sum_{m \in \mathcal{C}_k^\ell} \left| \frac{\mathrm{Cov}(f(V_m^{\ell+1}), g(V_m^{\ell+1}))}{\sigma(f(V_m^{\ell+1})) \cdot \sigma(g(V_m^{\ell+1}))} \right| \frac{|V_m^{\ell+1}|}{|V_k^\ell| \cdot \omega(V_m^{\ell+1})},
\end{equation}
where \( \mathrm{Cov} \) denotes the sample covariance and \( \sigma \) the sample standard deviation.
\end{itemize}
\end{definition}

\begin{remark} In \cite{coifman2013earth}, several metrics as in \eqref{eq:emd} are defined and proved equivalent
to Earth Mover's Distance with respect to the tree metric.
The parameters \( \alpha \), \( \beta \), and \( \epsilon \) modulate the effect of tree structure on the EMD. Increasing \( \alpha>0 \) places more weight on differences at coarser scales near the root, while \( \alpha = 0 \) yields uniform weighting across levels and \( \alpha < 0 \) emphasizes finer, leaf-level structure. The parameter \( \beta \) controls sensitivity to folder size: \( \beta > 0 \) favors larger folders, while \( \beta < 0 \) emphasizes smaller ones. The parameter \( \epsilon \) regulates the decay of the affinity; in practice, it is typically set as a constant multiple of the median EMD across all function pairs. A recent variant~\cite{lin2026joint} keeps this tree-EMD construction but replaces the parametric weight \( \omega(V_k^\ell) \) with one read from a hyperbolic embedding~\cite{nickel2017poincare} of the tree, so that the embedding geometry, rather than the fixed exponents \( \alpha \) and \( \beta \), sets the weight assigned to each folder.
\end{remark}

Now consider a kernel \( \mathcal{K} : Y \times X \to \mathbb{R} \). For convenience of notation, let \( W_{\mathcal{T}_X} \) and \( W_{\mathcal{T}_Y} \) denote tree-based affinities on  \( X \) and \( Y \), respectively. These affinities may be defined using either the multiscale Earth Mover's Distance \( W_{\mathcal{T}_X}^{\mathrm{EMD}} \), the multiscale correlation affinity \( W_{\mathcal{T}_X}^{\mathrm{corr}} \), as introduced in Definition~\ref{def_Tree_Aff}, or any other user-specified affinity measure constructed from the hierarchical structure of \( \mathcal{T}_X \) or \( \mathcal{T}_Y \).

\begin{definition}[Coupled Geometry]
Given a partition tree \( \mathcal{T}_X \) on \( X \), the \emph{dual affinity} between rows \( y_1, y_2 \in Y \) is defined as
\begin{equation}\label{eq:dual_row}
W_{\mathcal{T}_X}(y_1, y_2) := W_{\mathcal{T}_X}(\mathcal{K}_{y_1}, \mathcal{K}_{y_2}),
\end{equation}
where each \( \mathcal{K}_y(x) := \mathcal{K}(y, x) \) is viewed as a function on \( X \).

Similarly, given a partition tree \( \mathcal{T}_Y \) on \( Y \), the dual affinity between columns \( x_1, x_2 \in X \) is defined as
\begin{equation}\label{eq:dual_col}
W_{\mathcal{T}_Y}(x_1, x_2) := W_{\mathcal{T}_Y}(\mathcal{K}^{x_1}, \mathcal{K}^{x_2}),
\end{equation}
where \( \mathcal{K}^x(y) := \mathcal{K}(y, x) \) is viewed as a function on \( Y \).
\end{definition}

Using these constructions, the \texttt{Questionnaire} algorithm~\cite{ankenman2014geometry} refines the multiscale structure of the rows and columns and couples their geometries by alternating between the two axes, as detailed in Algorithm~\ref{alg:questionnaire_basic}.

\begin{algorithm}[H]
\caption{\texttt{Questionnaire}}\label{alg:questionnaire_basic}
\begin{algorithmic}[1]
\State \textbf{Input:} initial affinity \( W_X^{0} \) on \( X \); maximum number of iterations \( T \).
\For{\( t = 1, \dots, T \)}
    \State Build the column tree \( \mathcal{T}_X \) from \( W_X^{t-1} \).
    \State Set \( W_Y^{t} \gets W_{\mathcal{T}_X} \), the dual affinity \( \mathcal{T}_X \) induces on \( Y \), Eq.~\eqref{eq:dual_row}.
    \State Build the row tree \( \mathcal{T}_Y \) from \( W_Y^{t} \).
    \State Set \( W_X^{t} \gets W_{\mathcal{T}_Y} \), the dual affinity \( \mathcal{T}_Y \) induces on \( X \), Eq.~\eqref{eq:dual_col}.
\EndFor
\State \Return coupled trees \( \mathcal{T}_X, \mathcal{T}_Y \) and affinities \( W_X^{T}, W_Y^{T} \).
\end{algorithmic}
\end{algorithm}

This iterative refinement learns the latent hierarchical structure of the rows and columns, each tree being rebuilt from the dual affinity induced by the other. The tree-building steps of Algorithm~\ref{alg:questionnaire_basic} take an affinity matrix and return a multiscale partition tree; the \texttt{Questionnaire} framework does not prescribe how this is done, and any hierarchical tree construction may be used~\cite{ankenman2014geometry,coifman2011harmonic,mishne2017data}. In this work we build every tree with Algorithm~\ref{alg:affinity_ordering} (Section~\ref{sec:sfc}), which recursively bisects each folder by the Fiedler vector of the current affinity ($W_X^{t-1}$ for $\mathcal{T}_X$, $W_Y^{t}$ for $\mathcal{T}_Y$) and, in addition, fixes the left--right ordering of the sibling subtrees, and hence a one-dimensional order of the leaves (Section~\ref{sec:sfc}), automatically and as part of the same construction; bottom-up constructions such as the $k$-means-type tree of~\cite{coifman2011harmonic} leave this order arbitrary. We choose Fiedler bisection for its robust spectral properties and computational efficiency.

The framework has been developed in several directions. On the computational side, Ankenman~\cite{ankenman2014geometry} gives a practical implementation and experimental validation of the coupled geometry framework~\cite{pyquest}, and Coifman and Gavish~\cite{coifman2011harmonic} and Mishne et al.~\cite{mishne2017data} show that the resulting data-driven tree transforms yield highly sparse representations of the data. In applications, Georgiou et al.~\cite{georgiou2025clutter} use the \texttt{Questionnaire} algorithm to reconstruct disorganized and unlabeled datasets from chemical and bioengineering processes: on one-dimensional advection--diffusion systems and ensembles of Stuart--Landau oscillators, the iterative organization recovers meaningful spatial and parameter embeddings from scrambled data, enabling accurate data-driven identification and emulation of the system dynamics. On the theoretical side, this line of work explains why organizing a matrix to be smooth on its coupled trees makes it compressible. Coifman and Gavish~\cite[Thm.~9.1]{coifman2011harmonic} prove that if the expansion coefficients $c_i$ of a matrix in the tensor-Haar-like basis induced by its row and column trees have small $\ell_p$ entropy $e_p=\sum_i |c_i|^p$ ($0<p<2$), then the at most $\varepsilon^{-1}e_p$ coefficients that are large and supported on large rectangles approximate the matrix with an error controlled by $\varepsilon$. This notion of smoothness is defined by the trees alone: the tensor-Haar-like entropy does not depend on the left--right order of the folders, since exchanging two sibling folders only changes the sign of a coefficient. Smoothness in the sense of classical analysis, in which neighboring folders carry similar values, requires in addition an order of the folders. The left--right ordering of Algorithm~\ref{alg:affinity_ordering} supplies it, making the reorganized matrix a smooth function on the unit square (Section~\ref{sec:sfc}) that order-dependent transforms such as the \texttt{eGHWT} can compress (Figure~\ref{fig:q100_fixedpoints}). They also compute~\cite[Sec.~9.4.1]{coifman2011harmonic} that, for the dyadic partition tree of the torus $\mathbb{T}^d$ in its natural spatial coordinates, the tree-based correlation between the Fourier modes $e^{\mathrm{i}m\cdot x}$ and $e^{\mathrm{i}m'\cdot x}$ is, to second order, $1-c\,|m-m'|^2$, so that it recovers the Euclidean geometry of the frequency lattice $\mathbb{Z}^d$. This geometry exists only once the spatial tree is known. A random permutation of the rows and columns destroys that tree, and for an orthogonal matrix even the initial cosine affinity is the identity and carries no trace of $|m-m'|$ (Figure~\ref{fig:problem_dstiv}(c)). The coupled iteration nonetheless recovers it: in each iteration the rows improve their geometry using the tree on the columns, and the columns improve theirs using the tree on the rows, until the frequency geometry emerges even when the matrix is fully orthogonal (Figure~\ref{fig:sin_combined_vertical}). The natural coupled geometry thus acts as an attractor of the iteration, reached from random starts though not from every one (Figure~\ref{fig:q100_fixedpoints}), which is what the flip perturbation of Section~\ref{subsec:perturbation} addresses. Ankenman and Leeb~\cite{ankenman2018mixed} use wavelet shrinkage and Calder\'on--Zygmund-type decompositions to recover structured, piecewise-smooth components from matrices of mixed H\"older regularity, and prove that such matrices are efficiently approximated by thresholding tensor Haar-like basis functions supported on large rectangles. These guarantees concern tensor-Haar-like expansions, which suit Calder\'on--Zygmund-type kernels but not the oscillatory kernels that are the main targets of this work. For oscillatory kernels the relevant theory is that of the \texttt{Butterfly}, whose cost is set by the numerical ranks of the blocks: for the Fourier kernel and related special-function kernels, the numerical rank of a block to precision $\varepsilon$ is controlled by the frequency and by the area of the block in phase space~\cite{landau1962prolate,o2010algorithm}, so it stays bounded on the equal-area blocks of the factorization, and the interpolative decomposition controls the error of each block (\ref{app_ID} and~\ref{app_Butterfly}). These bounds hold in the natural coordinates of the operator; once the learned trees recover its phase-space geometry, they govern the compression of the reorganized matrix as well.

\section{Methodology}\label{sec:method}
Building on the reviewed \texttt{Questionnaire} and spectral-tree machinery, we introduce the components specific to this work: a spectral data-space-filling curve that resolves the tree-ordering ambiguity on which the downstream transforms depend, a stochastic \emph{flip perturbation} that keeps the iteration from stalling at a weak (poorly compressible) fixed point or periodic cycle, and the combined \texttt{QFFT} pipeline that couples geometry learning to compression.

\subsection{The Phase Space and the Spectral Data-Space-Filling Curve}\label{sec:sfc}

The idea that organizes this section is that of a phase space. The Fourier kernel $e^{-2\pi\mathrm{i}xy}$ is a function on the product of a space axis $x$ and a frequency axis $y$, the classical time--frequency, or phase, plane, and its fast algorithms exploit the geometry of that plane. By the uncertainty principle, a rectangle of the plane holding $N$ of the $N^2$ matrix entries carries only a bounded number of degrees of freedom~\cite{landau1962prolate}, so the kernel restricted to it has a numerical rank bounded independently of the position and shape of the rectangle; this is the property the \texttt{Butterfly} exploits~\cite{michielssen1996multilevel,o2010algorithm}, and the Walsh wavelet packets of~\ref{app_Walsh} are localized on tiles of the same plane. A general operator $\mathcal{K}$ also lives on a phase space, the product $Y\times X$ of its sensors (the range $Y$, indexing the rows) and its sources (the domain $X$, indexing the columns), in which the two are in duality, as positions and frequencies are for the Fourier transform: each is organized by the geometry that the matrix induces on the other. For the eigenvector matrices of Sections~\ref{sec:sh} and~\ref{sec:manifold} the sensors are eigenfunctions and the sources are points, so $Y\times X$ is again a space--frequency plane. Without coordinates, however, $Y\times X$ is only a set of index pairs. Our method learns its geometry from the matrix entries: the \texttt{Questionnaire} learns the two coupled trees, and Algorithm~\ref{alg:affinity_ordering} turns each tree into a coordinate.

This can be interpreted as the construction of a data-space-filling curve, a continuous, surjective mapping from the unit interval onto a higher-dimensional data-space, that adapts to the latent geometric structure of the data \cite{cantor1878beitrag}. Classical space-filling curves, such as the Hilbert~\cite{hilbert1935stetige}, Peano~\cite{peano1990courbe} and Z-order~\cite{morton1966computer} curves, traverse a regular grid by a fixed subdivision rule while approximately preserving proximity, and such locality-preserving traversals are used widely, from indexing structures and numerical simulation to cache-oblivious algorithms~\cite{bohm2020space}; but they require the coordinates of the grid. Our curve is a data-driven analog, built from the learned affinity alone, and Algorithm~\ref{alg:affinity_ordering} constructs it together with each tree of the \texttt{Questionnaire} iteration. The trees serve both compression schemes; the ordering serves the \texttt{eGHWT} and gives the phase space its analytic structure.

A tree becomes a coordinate through a left--right ordering. Definition~\ref{def:partition_tree} specifies an \emph{unordered} hierarchical grouping, which indices share a folder at each scale, but no one-dimensional order. A \emph{left--right ordering} of a binary tree $\mathcal{T}_X$ is a choice, at every internal folder, of which of its two children is the left one; it defines a total order on the points of $X$ recursively, each point of a left child preceding each point of its sibling. In the trees built by Algorithm~\ref{alg:affinity_ordering} the leaves are single points, so the order is on the points themselves; when a tree is truncated at leaf folders holding several points, as in the \texttt{Butterfly} (\ref{app_Butterfly}), the same rule orders the leaf folders, and the points inside each leaf keep the order of the full tree. A binary tree with $n$ internal folders admits $2^n$ left--right orderings; Algorithm~\ref{alg:affinity_ordering} selects one by the seam rule~\eqref{eq:seam} below, which fixes it up to a single global reflection.

Let $W$ be the current affinity on $X$; the construction on $Y$ is identical. For a folder $V\subseteq X$, \textsc{Order}$(V)$ returns an ordering of $V$, written as a list $\pi_V=(v_1,\dots,v_{|V|})$ of its points. It computes the Fiedler vector $\mathbf{f}^V$ of the subgraph $W(V,V)$, the minimizer of~\eqref{eq:fiedler} on $V$, and splits $V$ at the median $m_V$ of its entries into the two children
\begin{equation}\label{eq:split}
A=\{\,i\in V:\ f^V_i<m_V\,\},\qquad B=\{\,i\in V:\ f^V_i\ge m_V\,\}.
\end{equation}%
The children are ordered recursively, giving lists $\pi_A$ and $\pi_B$. For a list $\mathbf{u}=(u_1,\dots,u_p)$ let $\mathrm{id}(\mathbf{u})=\mathbf{u}$ and $\mathrm{rev}(\mathbf{u})=(u_p,\dots,u_1)$, and for a second list $\mathbf{w}=(w_1,\dots,w_q)$ let $\mathbf{u}\oplus\mathbf{w}=(u_1,\dots,u_p,w_1,\dots,w_q)$ be their concatenation. Each child's list is either kept or reversed, and the two are joined so that the points meeting at the seam are as close as possible in the Fiedler coordinate of the parent:
\begin{equation}\label{eq:seam}
\pi_V=\sigma_A(\pi_A)\oplus\sigma_B(\pi_B),\qquad
(\sigma_A,\sigma_B)=\argmin_{\sigma_A,\,\sigma_B\in\{\mathrm{id},\,\mathrm{rev}\}}\big|f^V_{a(\sigma_A)}-f^V_{b(\sigma_B)}\big|,
\end{equation}%
where $a(\sigma_A)$ is the last entry of $\sigma_A(\pi_A)$ and $b(\sigma_B)$ the first entry of $\sigma_B(\pi_B)$. Ties, which occur for instance when a child is a single point and $\mathrm{id}$ and $\mathrm{rev}$ coincide, are broken by a fixed order of the four choices, so that $\pi_V$ is well defined. The recursion stops at single points, so the leaves of the tree are singletons, and $\pi_X=\pi_V$ for $V=X$ is the left--right order of the tree.

\begin{algorithm}[!htb]
\caption{Construction of Tree and Space-Filling Curve}\label{alg:affinity_ordering}
\begin{algorithmic}[1]
\State \textbf{Input:} affinity matrix $W\in\mathbb{R}^{N\times N}$.
\State \textbf{Output:} balanced binary tree $\mathcal{T}$ and space-filling order $\pi$.
\Function{Order}{$V$}
    \If{$|V|\le 1$}
        \State \Return $V$
    \EndIf
    \State $\mathbf{f}\gets$ Fiedler vector of the normalized Laplacian of $W(V,V)$ (Eq.~\eqref{eq:fiedler} on $V$)
    \State $m\gets\operatorname{median}(\mathbf{f})$; \;\; $A\gets\{i\in V:\mathbf{f}_i<m\}$, \;\; $B\gets\{i\in V:\mathbf{f}_i\ge m\}$ (Eq.~\eqref{eq:split})
    \State $L\gets$ \Call{Order}{$A$}; \quad $R\gets$ \Call{Order}{$B$}
    \State reverse $L$ and/or $R$ so that $L$'s last point and $R$'s first point are closest in $\mathbf{f}$ (Eq.~\eqref{eq:seam})
    \State \Return $L$ followed by $R$
\EndFunction
\State $\pi\gets$ \Call{Order}{$\{1,\dots,N\}$}
\State $\mathcal{T}\gets$ the folders formed in the recursion; \quad \Return $\mathcal{T},\ \pi$
\end{algorithmic}
\end{algorithm}

The leaf order $\pi_X$ is a permutation of $X$, as in~\eqref{eq:problem}: $\pi_X(i)$ is the $i$-th point along the curve. Its inverse assigns to each point $x$ the curve coordinate $t(x):=(\pi_X^{-1}(x)-1)/N_X\in[0,1)$, and the coordinate $s(y)$ on $Y$ is defined in the same way from $\pi_Y$. Three properties of these coordinates follow from~\eqref{eq:split} and~\eqref{eq:seam}. First, the points of every folder are consecutive along the curve, and when $N_X=2^L$ the median split halves every folder, so the balanced tree maps each folder onto the grid points of a dyadic interval,
\begin{equation}\label{eq:dyadic}
t\big(V_k^\ell\big) = I_k^\ell\cap N_X^{-1}\mathbb{Z}, \qquad I_k^\ell := \big[k\,2^{-\ell},\,(k+1)\,2^{-\ell}\big), \qquad 0\le k<2^\ell,
\end{equation}%
where $k$ is the position of $V_k^\ell$ among the folders of level $\ell$, read from left to right; for other values of $N_X$ the intervals have nearly equal lengths. Two points of a common folder of level $\ell$ thus have curve coordinates within $2^{-\ell}$ of each other: closeness in the tree becomes closeness along the curve, as for the Hilbert curve, which traverses each dyadic subsquare of the square during one dyadic interval of its parameter~\cite{hilbert1935stetige}. Second, at every split the two halves are separated by the Fiedler coordinate of the parent, which minimizes the Dirichlet energy in~\eqref{eq:fiedler} among coordinates orthogonal to the constant; ordering by it is the spectral relaxation of placing high-affinity points next to each other on the line, that is, of minimizing the variation of the layout with respect to the learned affinity. Third, at every split~\eqref{eq:seam} chooses, among the four orientations of the two halves, the one that joins the two points closest in the Fiedler coordinate of the parent, so that neighbors along the curve tend to be neighbors in the learned affinity across folder boundaries as well as within folders. The seam rule is a local choice made once per folder; it does not minimize a global objective over all left--right orderings. Because the two axes are refined in alternation (Algorithm~\ref{alg:questionnaire_basic}), each iteration rebuilds the tree and its ordering on one axis from the geometry currently learned on the other, so that the row and column orderings improve together.

Through the two coordinates the reorganized matrix of~\eqref{eq:problem} becomes a function on the unit square, a coordinate model of the phase space $Y\times X$,
\begin{equation}\label{eq:phase}
\widetilde{\mathcal{K}}(s,t) := \mathcal{K}^{\pi}_{ji} \qquad\text{for}\quad s\in\Big[\frac{j-1}{N_Y},\frac{j}{N_Y}\Big),\quad t\in\Big[\frac{i-1}{N_X},\frac{i}{N_X}\Big),
\end{equation}%
and a pair of folders $V_{k'}^{\ell'}\subset Y$ and $V_k^{\ell}\subset X$ becomes the rectangle $I_{k'}^{\ell'}\times I_k^\ell$. The \texttt{Butterfly} tiles the square by such rectangles: at level $\ell$ it compresses the blocks $I_{k'}^{L-\ell}\times I_k^{\ell}$ (\ref{app_Butterfly}), all of the same area $2^{-L}$ (each holds $N$ entries when $N_X=N_Y=N=2^L$), as in the Fourier case; the problem statement~\eqref{eq:problem} asks that the reorganized kernel have low rank on these blocks. The \texttt{eGHWT} instead tiles the space--sequency plane of each axis: it expands $\widetilde{\mathcal{K}}$ in products of Walsh atoms, each localized on a unit-area tile of the space--sequency plane of its own axis and supported on a rectangle of the square, and selects the tiling of least cost (\ref{app_Walsh}); atoms of different sequency may share a support, so these rectangles overlap. The \texttt{Butterfly} blocks are products of folders and depend on the trees alone; the Walsh atoms depend also on the left--right order of the folders, which the curve fixes.

The curve is what makes this picture geometric rather than combinatorial. By the three properties above, a function that is nearly constant on the small folders of $\mathcal{T}_X$, the regularity measured by its tree-based Haar-like coefficients~\cite{coifman2011harmonic}, varies little along the curve. A kernel whose rows are regular on $\mathcal{T}_X$ and whose columns are regular on $\mathcal{T}_Y$, the bi-smoothness of~\cite{coifman2011harmonic} discussed in Section~\ref{sec:questionnaire}, becomes a regular function of both variables on the square, and its analytic structure becomes visible: the space--frequency plane of the DST-IV (Figure~\ref{fig:sin_combined_vertical}), the wavefronts of the Helmholtz kernel (Figures~\ref{fig:spiral_plan} and~\ref{fig:gh_matrices}), the dual index $(l,m)$ of the spherical harmonics (Figures~\ref{fig:SH_embedding} and~\ref{fig:SH_comparison}), and the phase space of the manifold eigenvectors (Figure~\ref{fig:Bean_comparison}). Since the reorganized matrix is then a function on the unit square, any transform of classical analysis can be carried out on the learned ordering. The \texttt{Butterfly} and the \texttt{eGHWT} are the two we use, and their compression of the reorganized matrices is the quantitative measure of the learned structure: the \texttt{Butterfly} measures the trees, and the \texttt{eGHWT} the trees together with their ordering (Figure~\ref{fig:q100_fixedpoints}). A fast transform calibrated to the discovered geometry, as the FFT is calibrated to the DFT, is a natural next step (Section~\ref{sec:discussion}).

\begin{figure}[!htb]
\centering
\includegraphics[width=0.85\linewidth]{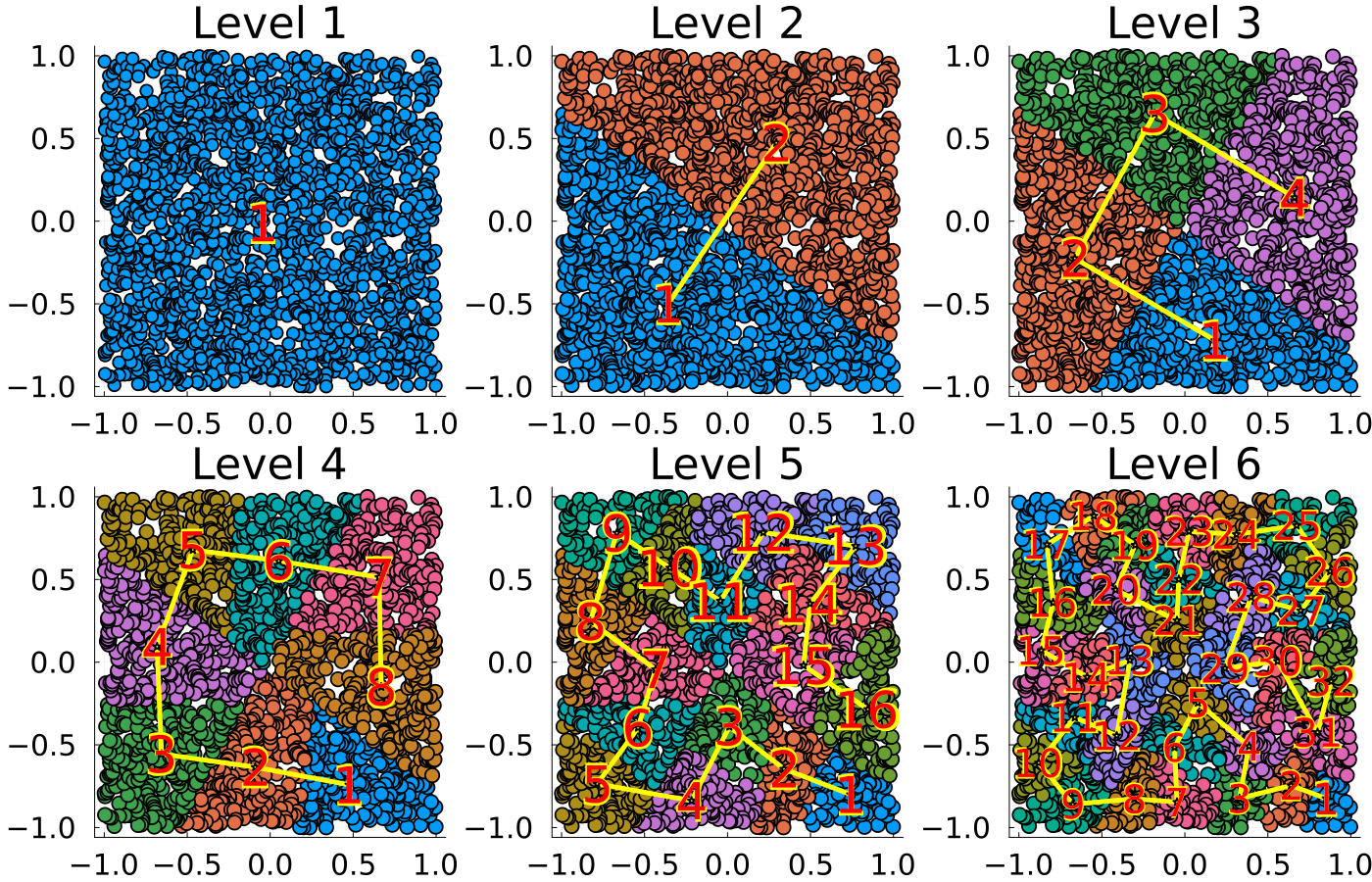}    \caption{Visualization of the space-filling curve at each tree level. Data points are uniformly sampled from the square \( [-1, 1] \times [-1, 1] \). Each subplot corresponds to a level of the partition tree. The title of each subplot indicates the tree level. The colors and the indices \( 1, 2, 3, \ldots \) denote the left-to-right ordering of the nodes at that level. The yellow line connects the nodes according to the space-filling curve traversal.}
    \label{fig:sfc_curve}
\end{figure}

In Figure~\ref{fig:sfc_curve}, we illustrate the performance of the space-filling curve on a synthetic dataset of \( 2048 \) points uniformly sampled on the square \( [-1, 1] \times [-1, 1] \). The affinity matrix \( W \) is computed using a Gaussian kernel based on pairwise Euclidean distances. As the level of the hierarchical partition increases, the space-filling curve progressively refines its path through the domain, revealing a smooth and coherent traversal that respects the spatial distribution of the points. At finer levels, the curve more densely fills the sampled surface, producing an ordering that captures the underlying geometry while maintaining spatial locality.

\subsection{Flip Perturbation}\label{subsec:perturbation}
\begin{figure}[!ht]
\centering
\includegraphics[width=\linewidth]{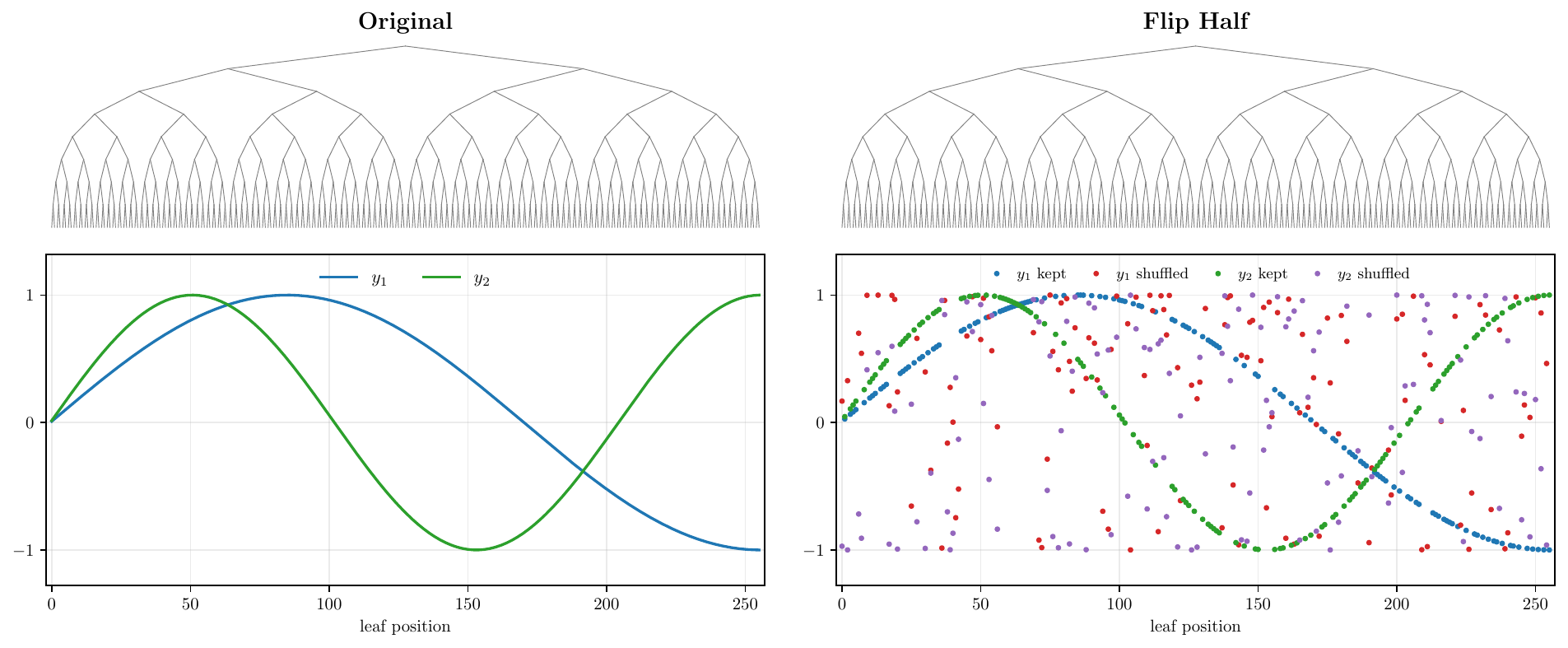}
\caption{Illustration of the flip perturbation ($N=256$). Two DST-IV rows, $y_1$ (row $y=1$) and $y_2$ (row $y=2$), lie along the leaves of a balanced binary tree. \textbf{Original}: the leaves are in order and both rows are smooth. \textbf{Flip Half}: a uniformly random half of the leaf positions is reshuffled and the same tree is rebuilt on the new order; the kept points stay on the smooth curves ($y_1$ blue, $y_2$ green), while the shuffled points scatter ($y_1$ red, $y_2$ purple). The kept half preserves the learned structure, while the shuffled half dislodges the current limit.}
\label{fig:flip}
\end{figure}

The \texttt{Questionnaire} iteration starts from an initial geometry, the seed affinity \( W_X^{0} \), and the limit it reaches depends on that geometry. As the demonstration in Figure~\ref{fig:q100_fixedpoints} shows, some initial affinities lead the iteration to a good, compressible fixed point, whereas others (for example, a noisy or otherwise uninformative initial affinity) drive it to a weak fixed point whose reorganization is far from smooth and hence poorly compressible, or to a short periodic cycle; either can be far from the most compressible organization the data admits. To escape such a  limit, we periodically apply a \emph{stochastic-flip perturbation} (Algorithm~\ref{alg:questionnaire}): every $\tau$ iterations, we permute a uniformly random half of the leaves of the tree built on each axis, holding the rest in place. For example, in Figure~\ref{fig:flip}, two smooth DST-IV rows lie along the leaves of a balanced binary tree; after a random half of the leaves are reshuffled, the kept points stay on the smooth curves while the shuffled points scatter. Because the two axes are refined in alternation, each perturbed tree induces a fresh dual affinity, and hence a fresh tree, on the other axis, so the flip propagates across the coupled geometry. The perturbed fraction is a free parameter; any ratio may be used, and here we use one half.

\begin{algorithm}[!htbp]
\caption{\texttt{Questionnaire} with flip perturbation}\label{alg:questionnaire}
\begin{algorithmic}[1]
\State \textbf{Input:} initial affinity \( W_X^{0} \) on \( X \); maximum number of iterations \( T\in\mathbb{N} \); perturbation period \( \tau\in\mathbb{N}\cup\{\infty\} \); compression precision \( \varepsilon\in(0,1) \).
\For{\( t = 1, \dots, T \)}
    \State Build the column tree \( \mathcal{T}_X \) and its leaf ordering from \( W_X^{t-1} \) by Algorithm~\ref{alg:affinity_ordering}.
    \If{\( t = n\tau \) for some integer \( n \)}
        \State Permute the points at a uniformly random half of the leaves of \( \mathcal{T}_X \), holding the rest in place.
    \EndIf
    \State Set \( W_Y^{t} \gets W_{\mathcal{T}_X} \), the dual affinity \( \mathcal{T}_X \) induces on the rows \( Y \), Eq.~\eqref{eq:dual_row}.
    \State Build the row tree \( \mathcal{T}_Y \) and its leaf ordering from \( W_Y^{t} \) by Algorithm~\ref{alg:affinity_ordering}.
    \If{\( t = n\tau \) for some integer \( n \)}
        \State Permute the points at a uniformly random half of the leaves of \( \mathcal{T}_Y \), holding the rest in place.
    \EndIf
    \State Set \( W_X^{t} \gets W_{\mathcal{T}_Y} \), the dual affinity \( \mathcal{T}_Y \) induces on the columns \( X \), Eq.~\eqref{eq:dual_col}.
\EndFor
\State \Return coupled trees \( \mathcal{T}_X, \mathcal{T}_Y \) and affinities \( W_X^{t^\star}, W_Y^{t^\star} \) from a chosen iterate \( t^\star \): the final iterate, or the one of least transform cost (the number of \texttt{eGHWT} coefficients or the \texttt{Butterfly} memory needed to represent the reordered matrix at precision \( \varepsilon \)).
\end{algorithmic}
\end{algorithm}

This perturbation works because, by the time the iteration reaches a limit, the learned organization already carries strong regularity with respect to the dual tree geometry: the rows and columns are coherently ordered and the reorganized matrix is smooth across most of its extent. Reshuffling only a fraction of the leaves keeps this coherence largely intact, so the perturbation does not by itself turn an inefficient limit into an efficient one; rather, it dislodges the current fixed point or periodic cycle while preserving the structure already discovered, and the iteration re-descends from a different basin toward a smoother, more compressible organization. The flip is therefore an optional step: when the initial affinity \( W_X^{0} \) is informative, the iteration usually descends directly to a compressible organization and the perturbation can be switched off (\( \tau=\infty \), which recovers the raw \texttt{Questionnaire} iteration).

\subsection{Questionnaire Factorization and Fast Transform}\label{section_kerenl_app}
Based on the previous sections,
our \texttt{QFFT} framework runs the \texttt{Questionnaire} with the flip perturbation (Algorithm~\ref{alg:questionnaire}) to learn coupled partition trees on the rows and columns of the input matrix. Within it, Algorithm~\ref{alg:affinity_ordering} builds each tree together with its leaf ordering from the current dual affinity, and reorganizing the rows and columns by the leaf orderings of the least-cost iterate reveals the joint structure of the matrix as an image. We then apply either the \texttt{Butterfly}, utilizing skeletonization built from row and column partition trees, or the \texttt{eGHWT} to the reordered matrix to obtain a compressed representation.
For the \texttt{Butterfly} algorithm, the compression is performed via interpolative decomposition (\texttt{ID}) maintaining a given precision \( \varepsilon \).
For \texttt{eGHWT}, compression is achieved by retaining the largest transform coefficients such that the total approximation error remains below a user-defined threshold \( \varepsilon \). In this work, we measure the approximation error using the Frobenius norm \( \|\cdot\|_F \) and $\ell_2$ norm \( \|\cdot\|_{\ell_2} \), as detailed in the following procedure.
The complete process is summarized in Algorithm~\ref{alg:compression}.

\begin{algorithm}[!t]
\caption{\texttt{QFFT}, Questionnaire Factorization and Fast Transform}
\label{alg:compression}
\begin{algorithmic}[1]
\Input 
Discrete matrix of kernel \( \mathcal{K} \in \mathbb{R}^{N_Y \times N_X} \); error parameter \( \varepsilon \in (0,1) \); initial affinity \( W_X^{0} \); maximum number of iterations \( T\in\mathbb{N} \); perturbation period \( \tau\in\mathbb{N}\cup\{\infty\} \); compression option: \texttt{"Butterfly"} or \texttt{"eGHWT"}

\State
Apply Algorithm~\ref{alg:questionnaire} with initial affinity $W_X^{0}$, iterations $T$, perturbation period $\tau$, and precision $\varepsilon$ to iteratively construct hierarchical partition trees $\mathcal{T}_X$ and $\mathcal{T}_Y$ and the affinities $W_X^{t^\star}, W_Y^{t^\star}$ of the returned iterate. Within this process, Algorithm~\ref{alg:affinity_ordering} is employed for tree construction; the final leaf orderings for $\mathcal{T}_X$ and $\mathcal{T}_Y$ are those it builds from $W_X^{t^\star}$ and $W_Y^{t^\star}$.
Reorganize \( \mathcal{K} \) according to the final row and column leaf orderings.

\If {compression option == \texttt{"Butterfly"}}
    \State Apply \texttt{Butterfly} factorization to the reordered matrix, skeletonizing the sub-blocks defined by the row and column trees $\mathcal{T}_X,\mathcal{T}_Y$, using \texttt{ID} with precision \( \varepsilon \) for every skeleton.
    \State \textbf{Output:} final-layer skeleton \( B^{(L)} \) and interpolation matrices \( \{P^{(\ell)}\} \).
\ElsIf {compression option == \texttt{"eGHWT"}}
    \State Apply the \texttt{eGHWT} on the row and column trees $\mathcal{T}_X,\mathcal{T}_Y$ to the reorganized \( \mathcal{K} \), and select the best tiling \( \mathcal{B}^* \) with basis \( \{\omega_{(q,p)}\}_{(q,p) \in \mathcal{B}^*} \) and coefficients \( c_{(q,p)} := \langle \mathcal{K}, \omega_{(q,p)} \rangle \).
    \State Sort \( |c_{(q,p)}| \) in descending order and retain the top \( n \) terms such that \( \sqrt{\sum_{i=n+1}^{N_Y N_X} c_{(q_i,p_i)}^2}\,\big/\,\norm{\mathcal{K}}_F \leq \varepsilon \).
    \State \textbf{Output:} compressed coefficients \( \{c_{(q_i,p_i)}\}_{i=1}^n \) and basis functions \( \{\omega_{(q_i,p_i)}\}_{i=1}^n \).
\EndIf
\end{algorithmic}
\end{algorithm}


We assume throughout this discussion that \( N_X, N_Y = O(N) \). In Algorithm~\ref{alg:compression}, if the compression option \texttt{Butterfly} is selected, the fast kernel application is simply
\begin{equation} 
 \widehat{\mathcal K}f(y_i) := \sum_{j_1} B^{(L)}_{i, j_{L-1}}
\sum_{j_2} P^{(L-1)}_{j_{L-1}, j_{L-2}}
\cdots
\sum_{j} P^{(0)}_{j_{0}, j} f(x_{j}).
\end{equation}
Provided that low-rank structure is successfully captured at each hierarchical level, the computational complexity of this procedure is \( \OO(N \log N) \).
If the compression option \texttt{eGHWT} is selected, the fast kernel application proceeds analogously through the \texttt{eGHWT} best basis: the input's \texttt{eGHWT} coefficients are aggregated against the retained kernel coefficients in the transform domain, and the output is reconstructed by the inverse \texttt{eGHWT}, again in \( \OO(N \log N) \) time when \( n = \OO(N\log N) \). We defer the derivation to \ref{app_Walsh}, equation~\eqref{eq_GHWT_apply}.

\section{Evaluation}\label{sec:eval}
{We evaluate the performance of our proposed \texttt{QFFT} framework on synthetic kernel matrices. Given a kernel matrix with rows and columns randomly permuted, we apply \texttt{QFFT} (Algorithm~\ref{alg:compression}), which learns the coupled geometry, reorganizes the rows and columns, and compresses the reorganized matrix with either \texttt{Butterfly} or the \texttt{eGHWT} at a fixed precision. 
We report the preparation time (\textit{Prep. Time}) of both the reorganization and compression stages, compare the memory of the original and compressed matrices, and assess the approximation accuracy and the computational efficiency of the resulting matrix-vector multiplication.} All entries of the matrices are stored in \texttt{Float64} format, corresponding to the standard 64-bit IEEE 754 double-precision floating-point representation. 
All code is implemented in \texttt{Python~3}. All experiments were conducted on a dual-socket system equipped with \texttt{AMD EPYC 7643} processors, providing a total of 96 physical CPU cores with a base clock frequency of 2.3~GHz.
To accelerate the \texttt{Questionnaire} preparation, which dominates the runtime, the computation of the tree-based affinity was parallelized across six {NVIDIA RTX A4500} GPUs, providing a total of 120~GB of GPU memory. The system used NVIDIA driver version 515.76 with CUDA~11.7 for hardware acceleration.
All other algorithms, including the \texttt{Butterfly} and \texttt{eGHWT}, were executed using single CPU thread.

Table~\ref{table:params} summarizes the six kernels studied below together with the algorithm parameters used for each: the maximum leaf-block size $\max_k|V_k^\ell|$ and interpolative-decomposition tolerance $\varepsilon$ for \texttt{Butterfly}, and the relative Frobenius tolerance $\varepsilon$ for \texttt{eGHWT}. The initial column affinity $W_X^0$ is a Gaussian kernel $W_X^0(i,i')=\exp(-\lVert z_i-z_{i'}\rVert^2/\sigma^2)$ on per-column feature vectors $z_i$, with $\sigma$ the median pairwise distance. When the kernel is informative the feature is the data itself, the $i$-th kernel column $z_i=\mathcal{K}^{x_i}$: this yields the \emph{Gaussian} affinity for the Green's and Helmholtz kernels, while the spherical-harmonic and manifold examples use the \emph{cosine} similarity $W_X^0(i,i')=\langle \mathcal{K}^{x_i},\mathcal{K}^{x_{i'}}\rangle/(\lVert \mathcal{K}^{x_i}\rVert\,\lVert \mathcal{K}^{x_{i'}}\rVert)$ of the same columns. For DST-IV, whose mutually orthogonal columns make $\mathcal{K}^{x_i}$ uninformative, the $z_i$ are instead drawn i.i.d.\ from a standard Gaussian (\emph{random}). For every example we apply Algorithm~\ref{alg:compression} to the original ordering, a random row/column permutation, and the reorganization learned by our method, and report preparation time, compressed memory, forward-transform time, and relative $\ell_2$ error. For every example we report results from a single seed, since $N$ grows large across these experiments. To avoid repetition, these parameter choices are collected here and not restated in each subsection.

\begin{table}[h]
\centering
\setlength{\tabcolsep}{4pt}
\resizebox{0.8\textwidth}{!}{
\begin{tabular}{l|cccccc}
\hline
 & \textbf{DST-IV} & \textbf{Green's} & \textbf{Helmholtz} & \textbf{Sph.\ (grid)} & \textbf{Sph.\ (rand.)} & \textbf{Manifold} \\
\hline
\textbf{Kernel} $\mathcal{K}$ & $\sin\!\big(\pi(y_j{-}\tfrac12)x_i\big)$ & $\dfrac{1}{\|x_i{-}y_j\|}$ & $\dfrac{\cos(2\pi\nu\|x_i{-}y_j\|)}{\|x_i{-}y_j\|}$ & $\mathcal{Y}_l^m$ & $\mathcal{Y}_l^m$ & eig.\ $\Delta_{\mathcal{M}}$ \\
\lightrule
\textbf{Tree aff.} & corr & EMD & EMD & corr & corr & corr \\
\lightrule
\textbf{Init.\ aff.} & random & Gaussian & Gaussian & cosine & cosine & cosine \\
\lightrule
$(T,\tau)$ & $(20,10)$ & $(10,\infty)$ & $(10,\infty)$ & $(10,\infty)$ & $(10,\infty)$ & $(10,\infty)$ \\
\lightrule
\texttt{Butterfly} & \multirow{2}{*}{$(128,10^{-11})$} & \multirow{2}{*}{$(64,10^{-11})$} & \multirow{2}{*}{$(64,10^{-11})$} & \multirow{2}{*}{$(64,10^{-11})$} & \multirow{2}{*}{$(64,10^{-11})$} & \multirow{2}{*}{$(256,10^{-6})$} \\
$(\max_k|V_k^\ell|,\varepsilon)$ & & & & & & \\
\lightrule
\texttt{eGHWT} $\varepsilon$ & $5{\times}10^{-2}$ & $10^{-3}$ & $10^{-3}$ & $5{\times}10^{-2}$ & $5{\times}10^{-2}$ & $5{\times}10^{-2}$ \\
\hline
\end{tabular}}
\caption{Summary of the numerical examples and the algorithm parameters used in each. ``Sph.\ (grid)'' and ``Sph.\ (rand.)'' are the spherical harmonics $\mathcal{Y}_l^m$ sampled on the structured grid and by random sampling, respectively, and the manifold kernel ``eig.\ $\Delta_{\mathcal{M}}$'' is the set of eigenvectors of the Laplace--Beltrami operator $\Delta_{\mathcal{M}}$ on the manifold. ``Tree aff.'' is the tree-based dual affinity on both the row and column trees $W_{\mathcal{T}_X},W_{\mathcal{T}_Y}$ driving \texttt{Questionnaire}, either the multiscale correlation (corr) or the tree-EMD (EMD) affinity; ``Init.\ aff.'' is the initial affinity $W_X^0$ seeding the iteration, defined in the text (random, Gaussian, or cosine); $(T,\tau)$ are the number of \texttt{Questionnaire} iterations and the perturbation period, with $\tau=\infty$ denoting no perturbation. For \texttt{Butterfly}, skeletonization begins at the finest tree level with at most $\max_k|V_k^\ell|$ columns per block, and each interpolative decomposition uses tolerance $\varepsilon$; for \texttt{eGHWT}, $\varepsilon$ is the relative Frobenius truncation tolerance.}
\label{table:params}
\end{table}

In the following tables, the parameter $N$ denotes the number of columns in the matrix to which our proposed algorithm is applied, while the number of rows is specified for each individual numerical example. Each table is partitioned into distinct blocks for the \texttt{Butterfly}, \texttt{eGHWT}, and \texttt{Questionnaire} frameworks. The ``Prep. Time'' row reports the computation time of each algorithm across various matrix sizes in seconds (sec). For both the \texttt{Butterfly} and \texttt{eGHWT} methods, the subsequent rows (\textit{Permuted}, \textit{Original}, and \textit{Reorganized}) report the memory consumption in MiB after applying the factorization under different ordering schemes: a random permutation of rows and columns, the original ordering of the row and column spaces, and the matrix reorganized by our proposed \texttt{QFFT} method. The efficiency of the kernel application is evaluated through two metrics: ``Direct Evaluation'' lists the time required to perform a standard dense matrix--vector multiplication in seconds, while ``Forward Transform'' reports the time in seconds required to apply the factorized kernel using the accelerated algorithms described in Section~\ref{section_kerenl_app}. Finally, the ``Relative $\ell_2$ Error'' column quantifies the numerical accuracy of the application. Given a test function $f \in \mathbb{R}^N$ and the approximate image $\widehat{\mathcal{K}}f$, obtained via the \texttt{Butterfly} method~\eqref{eq_bf_apply} or the \texttt{eGHWT} method~\eqref{eq_GHWT_apply}, the error is defined as 
\begin{equation}
    \frac{\|\mathcal{K}f - \widehat{\mathcal{K}}f\|_{\ell_2}}{\|\mathcal{K}f\|_{\ell_2}},
\end{equation}
computed using test functions drawn uniformly at random from the interval $[0, 1]$ and averaged over 300 trials.  Since the best basis selection algorithm of eGHWT requires $\OO(N^2 (\log N)^2)$ memory, cases where this exceeds our CPU memory capacity are indicated by ``--'' in the table.

The two compressors target different regimes, which fixes the tolerances in Table~\ref{table:params}. \texttt{Butterfly}, developed for oscillatory kernels, attains high-precision approximation with both memory and application cost $\OO(N\log N)$, and runs at a tight tolerance throughout ($\varepsilon=10^{-11}$, and $10^{-6}$ for the manifold). \texttt{eGHWT} instead builds localized, data-adaptive wavelet-like bases on the learned row and column trees and represents the kernel sparsely in that best basis; its application is likewise $\OO(N\log N)$, but its piecewise-constant $\pm1$ bases suit smooth or piecewise-structured operators and are poorly matched to oscillatory kernels. Its tolerance is therefore set per example: $10^{-3}$ for the milder Green's and Helmholtz kernels and a looser $5\times10^{-2}$ for the strongly oscillatory DST-IV, spherical-harmonic, and manifold examples, whose \texttt{eGHWT} memory grows toward $\OO(N^2)$ and is omitted (the ``--'' entries) at the largest $N$. In either case $\varepsilon$ is a genuine accuracy control, but it trades accuracy chiefly against \emph{memory} (and the $\OO(N\log N)$ application), not preparation time. The achieved relative $\ell_2$ error tracks $\varepsilon$ (reported for every example in the tables below). A tighter $\varepsilon$ raises the \texttt{Butterfly} retained per-block rank and the \texttt{eGHWT} retained coefficient count, and hence the compressed memory and the application cost; but the preparation is dominated by $\varepsilon$-independent work, the \texttt{Questionnaire} iteration ($\OO(N^3)$) and the compression step itself (the \texttt{eGHWT} best-basis selection, $\OO(N^2(\log N)^2)$, or the \texttt{Butterfly} construction, $\OO(N^2\log N)$), so relaxing $\varepsilon$ lowers the memory appreciably while reducing the preparation time only modestly.

Finally, we state the break-even explicitly. The preparation is a large \emph{one-time} cost, running to thousands of times a single dense matrix--vector product; moreover, the \texttt{Questionnaire} preparation is GPU-accelerated whereas the direct evaluation is single-threaded, so on identical hardware the preparation would be costlier still (Table~\ref{tab:timing} reports single-thread timings for a like-for-like measurement). The method is therefore not meant to outrun an existing fast transform on a single application; its purpose is to \emph{discover} an $\OO(N\log N)$ transform for an operator that has none in its given ordering. Once prepared, each application costs $\OO(N\log N)$ rather than the $\OO(N^2)$ of direct evaluation, so the one-time cost is amortized only when the operator is applied many times, or at large $N$ where the per-application gap is decisive; the break-even, in number of applications, is correspondingly large.

\section{Numerical Simulations}

For each example below we introduce the kernel and its construction, describe the geometry the \texttt{Questionnaire} learns, and report the compression in a table of preparation time, memory, forward-transform time, and error, with a plot of time and memory against $N$ beneath it. As in the demonstration of Section~\ref{sec:dstiv}, the one-time preparation dominates the cost: each \texttt{Questionnaire} iteration scales as $\OO(N^3)$ and the \texttt{eGHWT} best-basis selection as $\OO(N^2(\log N)^2)$, whereas once built the operator applies to each vector in only $\OO(N\log N)$, below the $\OO(N^2)$ dense multiplication; our single-threaded \texttt{eGHWT} forward transform, however, remains slower than dense multiplication in absolute terms. Across every example the reorganization discovers the geometric structure hidden in the scrambled matrix, compressing far better than a random permutation, and the order it finds need not be a known natural one. Where a fast $\OO(N\log N)$ structure is present the iteration recovers it; for noisier data, such as the randomly sampled spherical harmonics and the manifold Laplacian, the matrix may not be compressible all the way to $\OO(N\log N)$, yet the learned order still yields a transform faster than dense multiplication once $N$ is large enough.

\subsection{Discrete Sine Transform (DST-IV)}\label{sec_sine}

This example is the DST-IV basis $\mathcal{K}_{ji}=\sin\!\big(\pi(y_j-\tfrac12)x_i\big)$. Earlier demonstrations of the \texttt{Questionnaire} start from a known, informative geometry: Coifman and Gavish~\cite[Sec.~9.4.2]{coifman2011harmonic}, for instance, organize $200$ sines and cosines sampled at $512$ points, an oversampled, non-square matrix whose columns are not orthogonal and therefore carry an informative initial affinity. Here the DST-IV is square and orthogonal in both its rows and columns, so the cosine affinity is the identity on both axes and carries no geometry; we therefore seed the iteration with the \emph{random} Gaussian affinity of Section~\ref{sec:eval}. As its demonstration in Section~\ref{sec:dstiv} shows, the \texttt{Questionnaire} iteration with the flip perturbation (Algorithm~\ref{alg:questionnaire}) nonetheless recovers the one-dimensional structure of the space $X$ and the frequency $Y$ (Figure~\ref{fig:sin_combined_vertical}). Table~\ref{table:sinkx2} reports the compression results across sizes $N$.

\begin{table}[h]
\centering
\makebox[\textwidth][c]{%
\begin{minipage}[c]{0.52\textwidth}
\centering
{\scriptsize\renewcommand{\arraystretch}{0.85}\setlength{\tabcolsep}{2pt}%
\resizebox{\linewidth}{!}{%
\begin{tabular}{lcccccc}
    \hline
    \textbf{N} & \textbf{1024} & \textbf{2048} & \textbf{4096} & \textbf{8192} & \textbf{16384} & \textbf{32768} \\
    Direct Evaluation (sec) & 5.71E-4 & 6.86E-3 & 1.66E-2 & 6.80E-2 & 3.30E-1& 1.29E0 \\
    \hline
       {\texttt{Questionnaire}}\\
    Prep. Time (sec)     & 8.90E1  & 1.15E2  & 2.76E2 & 1.05E3  & 4.47E3  & 1.30E4  \\
    \hline
{\texttt{Butterfly}} \\
    Permuted (MiB)      & 1.80E1     & 8.40E1     & 2.96E2    & 1.36E3   & 4.77E3   & 1.89E4
    \\
    Original (MiB)      & 9.19E0    & 2.47E1   & 6.29E1   & 1.54E2  & 3.65E2  & 8.51E2\\
    Reorganized (MiB)   & 6.59E0    & 2.06E1   & 6.32E1   & 1.50E2    & 4.26E2  & 1.45E3     \\
    Prep. Time (sec)      & 1.40E0  & 9.24E0  & 1.89E2  & 5.70E2  & 1.53E3  & 4.03E3 \\
    Forward Transform (sec)      & 1.05E-3  & 2.45E-3  & 8.60E-3  & 1.82E-2  & 3.71E-2  & 6.24E-2  \\
    Relative $\ell_2$ Error  & 9.42E-12  & 2.09 E-11 &  1.54E-11  & 1.45E-11 & 1.12E-11  & 1.75E-11  \\ \hline
    {\texttt{eGHWT}}\\
    Permuted (MiB)      & 1.43E1     & 5.67E1     & 2.25E2    & 8.98E2  & 3.59E3   & --     \\
    Original (MiB)      & 2.71E0    & 8.45E0   & 2.66E1   & 8.48E1  & 2.65E2  & --    \\ Reorganized (MiB)     & 1.00E1   & 3.83E1   & 1.58E2    & 6.38E2  & 1.98E3 &   --  \\
    Prep. Time (sec)      & 2.98E1  & 1.77E2  & 7.30E2  & 4.31E3  & 2.28E4 & --  \\
     Forward Transform (sec)    & 1.76E-1  & 9.69E-1  & 3.61E0  & 1.40E1  & 4.33E1 & --  \\
    Relative $\ell_2$ Error  & 3.92E-2  & 4.92E-2  & 4.92E-2  & 4.97E-2  & 4.70E-2  & -- \\ \hline
    \end{tabular}}}
\end{minipage}\hspace{0.045\textwidth}
\begin{minipage}[c]{0.52\textwidth}
\centering
\includegraphics[width=\linewidth]{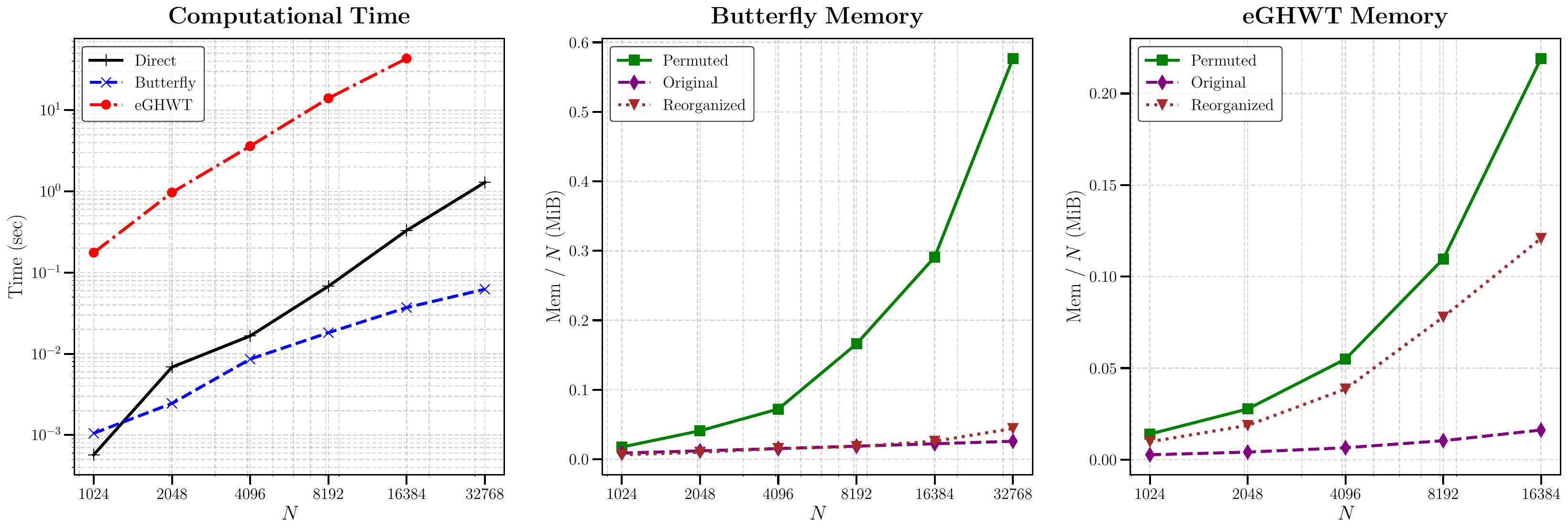}
\end{minipage}}
\caption{DST-IV: times, errors, and memory usage (left), and computational time and memory versus $N$ (right).}
\label{table:sinkx2}
\end{table}

\subsection{Acoustic Waves and Potential Operator}
\begin{figure}
    \centering
\includegraphics[width=0.95\linewidth,trim=3cm 1cm 1cm 0cm, clip]{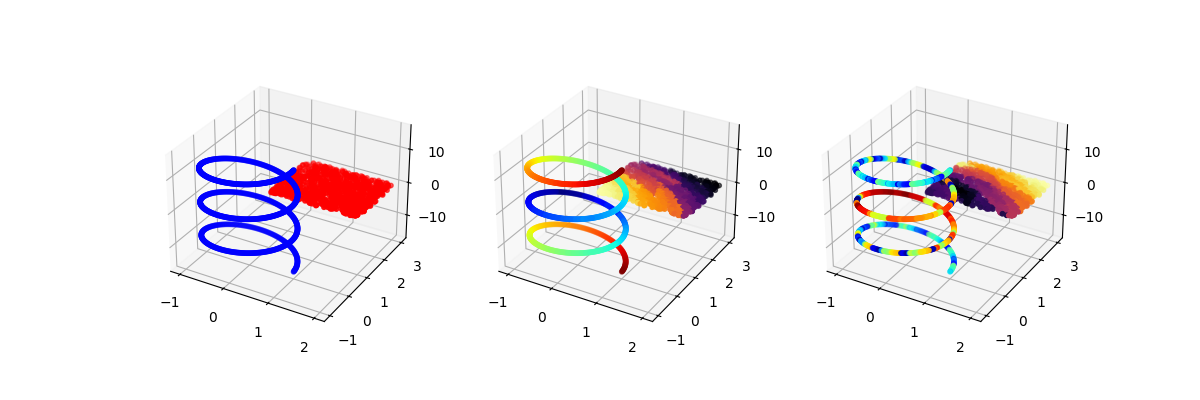}
    \caption{
Left: spatial layout of the helix (points \( \{x_i\}_{i=1}^N \)) and the flat sheet (points \( \{y_j\}_{j=1}^N \)). 
Middle: same layout, with points colored according to tree-based space filling curve learned via Algorithm~\ref{alg:affinity_ordering} for the kernel $\frac{1}{\norm{x_i-y_j}}$.
Right: same layout, with points colored according to tree-based space filling curve learned via Algorithm~\ref{alg:affinity_ordering} for the kernel $\frac{\cos(2\pi \norm{x_i-y_j})}{\norm{x_i-y_j}}$.}
\label{fig:spiral_plan}
\end{figure}

Consider an acoustic wave propagation scenario between two point clouds. We construct the column and row spaces embedded in \( \mathbb{R}^3 \): a one-dimensional helix \(X =  \{x_i\}_{i=1}^N \) (blue, the columns) and a two-dimensional sheet \( Y = \{y_j\}_{j=1}^N \) (red, the rows), as illustrated in the left panel of Figure~\ref{fig:spiral_plan}. The "acoustic" interaction between source and target points is modeled by the kernel
\begin{equation}
\mathcal K_{ji} = \frac{\cos(2\pi \nu \|x_i - y_j\|)}{\|x_i - y_j\|},
\end{equation}
where \( \nu \) is the frequency parameter. This kernel corresponds to the real part of the three-dimensional Helmholtz equation with wavenumber \( k = 2\pi \nu \). It captures oscillatory wave behavior with spatial decay and is commonly used to model time-harmonic acoustic interactions in free space as a potential operator.
When \( \nu = 0 \), the kernel simplifies to
\begin{equation}
\mathcal K_{ji} = \frac{1}{\|x_i - y_j\|},
\end{equation}
which corresponds to the static Green’s function of the Laplace equation in three dimensions. In this case, the kernel models a classical potential operator describing steady-state interactions such as electrostatic or gravitational forces, with long-range influence and no oscillatory behavior. We examine the scenario with the equispaced helical source points \( \{x_i\} \) and the non-equispaced target points \( \{y_j\} \) on the sheet, both are randomly shuffled and remain unorganized. The static case of this helix--sheet layout is the potential-operator example of Coifman and Gavish~\cite[Secs.~9.3 and~9.8]{coifman2011harmonic}, who recover the coupled geometry of the scrambled $600\times600$ matrix and compress it in a tensor-Haar-like basis. Such a piecewise-constant basis suits Calder\'on--Zygmund-type potentials, whose interactions are smooth away from the near field, but not oscillatory kernels. We extend the example to the oscillatory Helmholtz kernel, to fast matrix--vector multiplication by the \texttt{Butterfly} and \texttt{eGHWT}, and to sizes up to $N=32768$.

\begin{figure}
    \centering
    \makebox[\textwidth][c]{%
    \begin{minipage}[b]{0.52\textwidth}
        \centering
        \includegraphics[width=\linewidth]{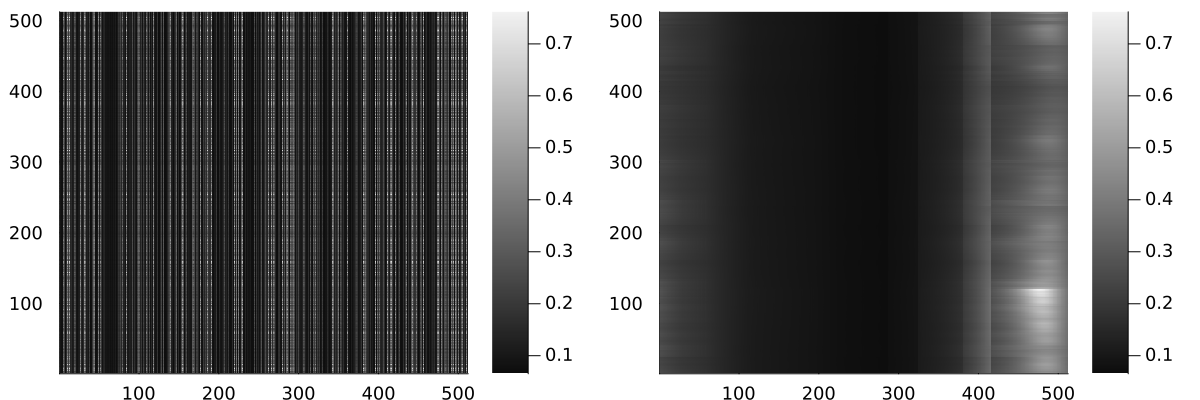}
        \subcaption*{Green function $\frac{1}{\norm{x-y}}$}
    \end{minipage}\hspace{0.045\textwidth}
    \begin{minipage}[b]{0.52\textwidth}
        \centering
        \includegraphics[width=\linewidth]{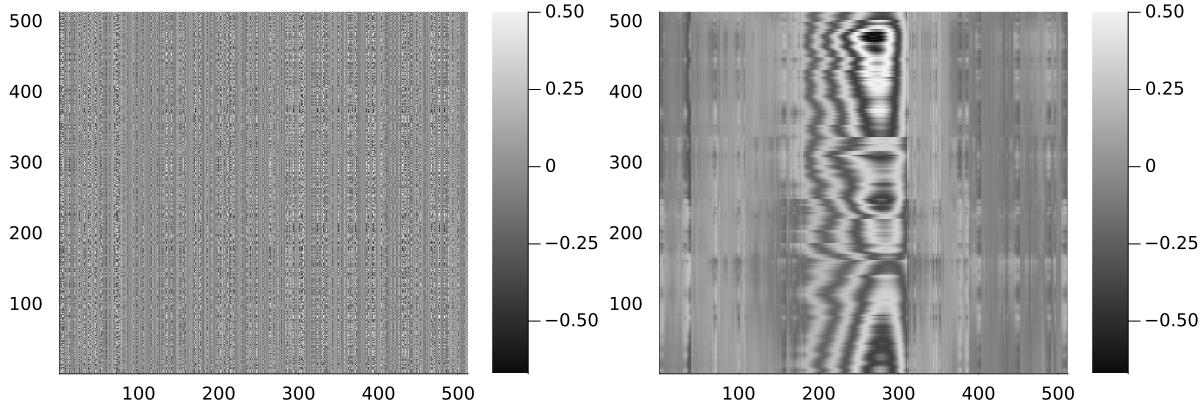}
        \subcaption*{Helmholtz kernel $\frac{\cos(2\pi\norm{x-y})}{\norm{x-y}}$}
    \end{minipage}}
    \caption{Unorganized matrix (left of each panel) and after reorganizing the rows and columns (right), for the Green function and the Helmholtz kernel.}
    \label{fig:gh_matrices}
\end{figure}

We recover the source--target correspondence directly from the permuted matrix. The middle and right panels of Figure~\ref{fig:spiral_plan} color the helix and sheet by the learned curve for the Green ($\nu=0$) and Helmholtz ($\nu=1$) kernels. Because the Helmholtz kernel oscillates in the source--target distance $\|x_i-y_j\|$, the learned curve arranges the helix and sheet so that this distance, and hence the wave's phase, varies coherently along the ordering (right panel): the oscillation that appears as noise in the permuted matrix is reorganized into smooth wavefronts, a discrete ``geometric optics'' read directly off the kernel. Figure~\ref{fig:gh_matrices} shows the potential-operator and acoustic-wave matrices before and after reorganization.

Table~\ref{table:gh_results} reports the compression results across sizes $N$.

\begin{table}[t]
\centering
\makebox[\textwidth][c]{%
\begin{minipage}[t]{0.52\textwidth}
\centering
\subcaption*{Green function $\frac{1}{\norm{x-y}}$}
{\scriptsize\renewcommand{\arraystretch}{0.85}\setlength{\tabcolsep}{2pt}%
\resizebox{\linewidth}{!}{%
\begin{tabular}{lcccccc}
    \hline
    \textbf{N} & \textbf{1024} & \textbf{2048} & \textbf{4096} & \textbf{8192} & \textbf{16384} & \textbf{32768} \\
    Direct Evaluation& 7.12E-4 & 3.68E-3 & 1.78E-2 & 8.10E-2 & 3.31E-1 &1.32E0\\
    \hline
       {\texttt{Questionnaire}}\\
    Prep. Time (sec)     & 5.58E1  & 9.16E1  & 2.28E2 & 1.00E3  & 2.51E3  & 1.20E4\\
    \hline
{\texttt{Butterfly}} \\
    Permuted (MiB)      & 1.93E0 & 4.52E0 & 1.03E1 & 2.29E1 & 5.07E1 & 8.46E1
    \\
    Reorganized (MiB)   & 5.34E-1 & 9.89E-1 & 1.83E0 & 3.38E0 & 6.22E0 & 1.14E1     \\
    
    Prep. Time (sec)      & 4.71E-1  & 7.35E-1  & 3.26E0  & 1.55E1  & 6.31E1  & 3.38E2 \\
    Forward Transform (sec)      & 3.53E-4  & 1.16E-3  & 2.60E-3  & 5.91E-3  & 1.21E-2  & 3.48E-2  \\
    
    Relative $\ell_2$ Error  & 3.81E-12  & 
    3.44E-12 &  3.68E-12  & 4.05E-12 & 4.20E-12  &  4.23E-12 \\ \hline
    {\texttt{eGHWT}}\\
    Permuted (MiB)      & 1.31E1     & 5.12E1     & 2.04E2    & 8.18E2    & 3.30E3  & --    \\
    Reorganized (MiB)     & 4.72E-1   & 7.41E-1   & 1.03E0    & 1.38E0  & 1.83E0 &  --   \\ 
    Prep. Time (sec)      & 4.45E1  & 2.00E2  & 7.63E2  & 4.04E3  & 2.27E4 & --  \\
     Forward Transform (sec)    & 1.01E-1  & 2.12E-1  & 4.44E-1  & 9.42E-1 & 1.99E0 &--  \\
    Relative $\ell_2$ Error  & 1.15E-4  & 1.22E-4  & 1.31E-4  & 1.48E-4  & 1.68E-4  & --  \\ \hline
 
    \end{tabular}}}

\vspace{4pt}
\includegraphics[width=\linewidth]{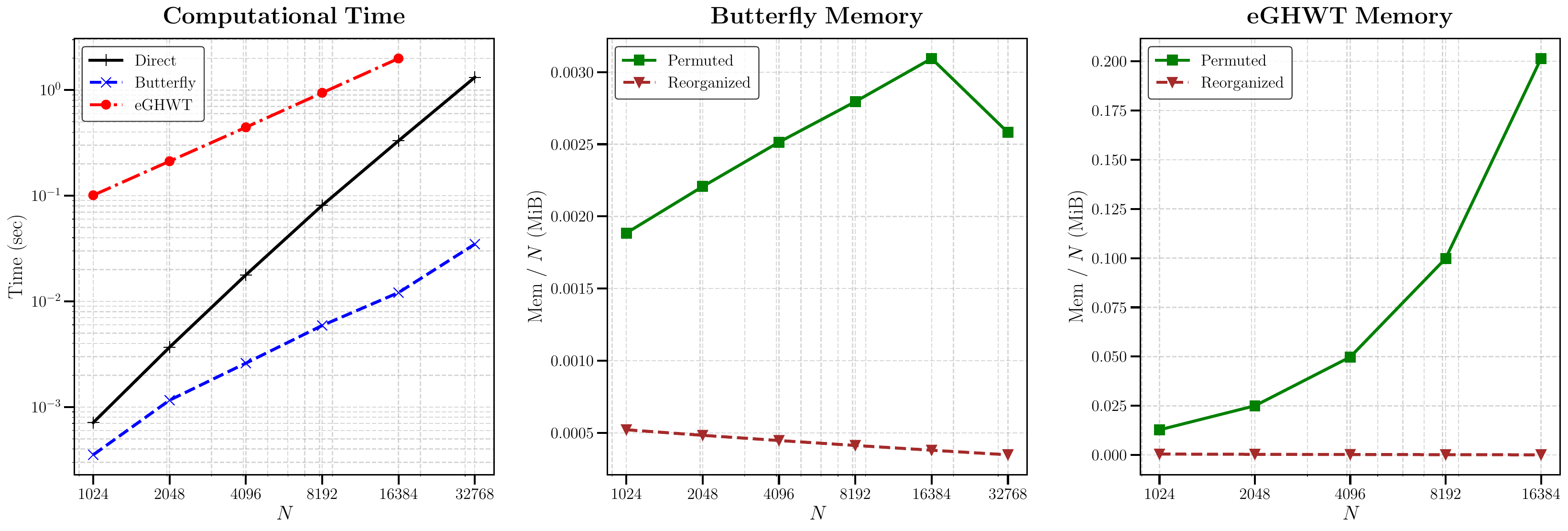}
\end{minipage}\hspace{0.045\textwidth}
\begin{minipage}[t]{0.52\textwidth}
\centering
\subcaption*{Helmholtz kernel $\frac{\cos(2\pi\norm{x-y})}{\norm{x-y}}$}
{\scriptsize\renewcommand{\arraystretch}{0.85}\setlength{\tabcolsep}{2pt}%
\resizebox{\linewidth}{!}{%
\begin{tabular}{lcccccc}
    \hline
    \textbf{N} & \textbf{1024} & \textbf{2048} & \textbf{4096} & \textbf{8192} & \textbf{16384} & \textbf{32768} \\
    Direct Evaluation & 9.14E-4 & 6.08E-3 & 2.56E-2 & 8.13E-2 & 3.33E-1 & 1.32E0\\
    \hline
       {\texttt{Questionnaire}}\\
    Prep. Time (sec)     & 5.61E1  & 8.78E1  & 2.74E2 & 1.14E3  & 2.67E3  & 1.60E4   \\
    \hline
{\texttt{Butterfly}} \\
    Permuted (MiB)      & 8.15E0 & 2.09E1 & 5.03E1 & 1.18E2 & 2.69E2 & 6.04E2
    \\
    Reorganized (MiB)   & 3.50E0 & 7.94E0 & 1.25E1 & 2.41E1 & 4.02E1 & 6.27E1     \\
    
    Prep. Time (sec)      & 4.85E-1  & 1.55E0  & 4.94E0  & 2.22E1  & 7.68E1  &  3.01E2\\
    Forward Transform (sec)      & 4.61E-4  & 1.96E-3  & 3.78E-3  & 8.02E-3  & 1.62E-2  & 3.79E-2  \\
    
    Relative $\ell_2$ Error  & 2.43E-11  & 
    2.67E-11 &  2.27E-11  & 2.71E-11 & 2.66E-11  & 2.28E-11  \\ \hline
    {\texttt{eGHWT}}\\
    Permuted (MiB)      & 2.22E1     & 8.91E1     & 3.56E2    & 1.42E3  & 5.70E3   & --     \\
    Reorganized (MiB)     & 1.32E1   & 3.81E1   & 9.40E1    & 2.19E2  & 4.48E2 &  --   \\ 
    Prep. Time (sec)      & 3.18E1  & 1.72E2  & 7.38E2  & 3.79E3  & 2.27E4 & --  \\
     Forward Transform (sec)    & 3.53E-1  & 9.81E-1  & 2.36E0  & 5.41E0  & 1.14E1 & --  \\
    Relative $\ell_2$ Error  & 2.89E-4  & 2.07E-4  & 1.61E-4  & 1.37E-4 & 1.2E-4  & --  \\ \hline
 
    \end{tabular}}}

\vspace{4pt}
\includegraphics[width=\linewidth]{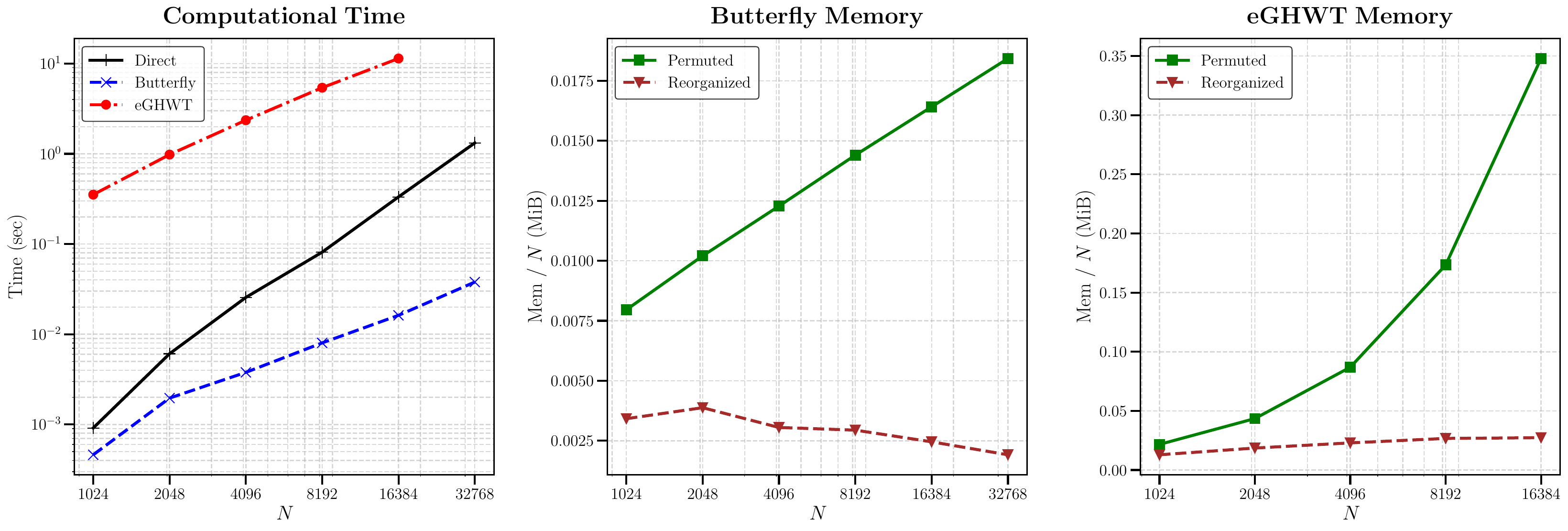}
\end{minipage}}
\caption{Times, errors, and memory (tables) and computational time and memory versus $N$ (plots) for the Green function (left) and the Helmholtz kernel (right).}
\label{table:gh_results}
\end{table}

\subsection{Spherical Harmonics}\label{sec:sh}
We organize the spherical-harmonic eigenfunctions $\mathcal{Y}_l^m$ of the Laplace--Beltrami operator on $\mathbb{S}^2$, indexed by degree $l\in\mathbb{N}\cup\{0\}$ and order $-l\le m\le l$,
\begin{equation}
\Delta_{\mathbb{S}^2} \mathcal{Y}_l^m = -l(l+1) \mathcal{Y}_l^m,
\end{equation}
with $(l_{\max}+1)^2$ harmonics up to degree $l_{\max}$. In real form, with colatitude $\theta\in[0,\pi]$ and azimuth $\phi\in[0,2\pi)$,
\begin{equation}
\mathcal{Y}_l^m(\theta, \phi) =
\begin{cases}
\sqrt{2}\, N_{l}^{|m|}\, \mathcal{P}_l^{|m|}(\cos\theta) \cos(m\phi), & m > 0, \\
N_{l}^0\, \mathcal{P}_l^0(\cos\theta), & m = 0, \\
\sqrt{2}\, N_{l}^{|m|}\, \mathcal{P}_l^{|m|}(\cos\theta) \sin(|m|\phi), & m < 0,
\end{cases}
\end{equation}
where $\mathcal{P}_l^m$ are the associated Legendre polynomials and $N_{l}^m=\sqrt{\tfrac{2l+1}{4\pi}\tfrac{(l-m)!}{(l+m)!}}$.

The eigenfunctions form a basis matrix $\mathcal{K}\in\mathbb{R}^{(N/4)\times N}$ whose row space is $\bigoplus_{l=0}^{l_{\max}}\mathcal{H}_l$, with $\mathcal{H}_l=\mathrm{span}\{\mathcal{Y}_l^m\}_{m=-l}^{l}$ the degree-$l$ eigenspace. We use two sampling schemes that differ both in the points and in the basis chosen within each $\mathcal{H}_l$. (1)~\emph{Structured grid}: the canonical harmonics $\mathcal{Y}_l^m$ evaluated on the equiangular lattice $\theta_a=(a+\tfrac12)\pi/N_\theta$, $\phi_b=2\pi b/N_\phi$ with $N_\theta=N_\phi=2(l_{\max}+1)$, giving $N=4(l_{\max}+1)^2$ points. (2)~\emph{Random sampling}: $N$ points $x_i=g_i/\lVert g_i\rVert_2$ with $g_i\sim\mathcal{N}(0,I_3)$ drawn independently, which are uniform on the sphere; here $l_{\max}=\lfloor-1+\sqrt{N/4}\rfloor$ gives about $N/4$ eigenfunctions. Because the $2l+1$ harmonics of a degree span a single eigenspace $\mathcal{H}_l$ (all sharing the eigenvalue $-l(l+1)$), every orthonormal basis of $\mathcal{H}_l$ is an equally valid set of degree-$l$ eigenfunctions; for the random scheme we therefore take a Haar-random orthonormal basis of each $\mathcal{H}_l$, so that no azimuthal order is privileged. A fast spherical-harmonic transform presupposes a structured quadrature grid such as the one above; it does not apply to the randomly sampled points, which are not a permutation of any such grid, so the randomly sampled harmonics have no dedicated fast transform to fall back on.

For the structured grid we order the rows by azimuthal order $m$, and within each $m$ by increasing degree $l$, using the triangular number $S(n)=\tfrac{n(n+1)}{2}$,
\begin{equation}
    j(l, m) =
    \begin{cases}
    S(l_{\max} + m) + l + m, & m < 0, \\
    S(l_{\max}) + S(l_{\max} + 1) - S(l_{\max} - m + 1) + l - m, & m \ge 0,
    \end{cases}
\end{equation}
which clusters equal-frequency harmonics and exposes the block-diagonal structure under rotational symmetry; for the randomly rotated scheme the order $m$ is undefined, so the rows are grouped only by degree $l$. For the grid the column index is $i(a,b)=aN_\phi+b$ (longitudinal scans at fixed latitude); for random sampling, columns follow the sampling order.

\begin{figure}
    \centering
    \begin{subfigure}[b]{\linewidth}
         \centering
         \includegraphics[width=\linewidth]{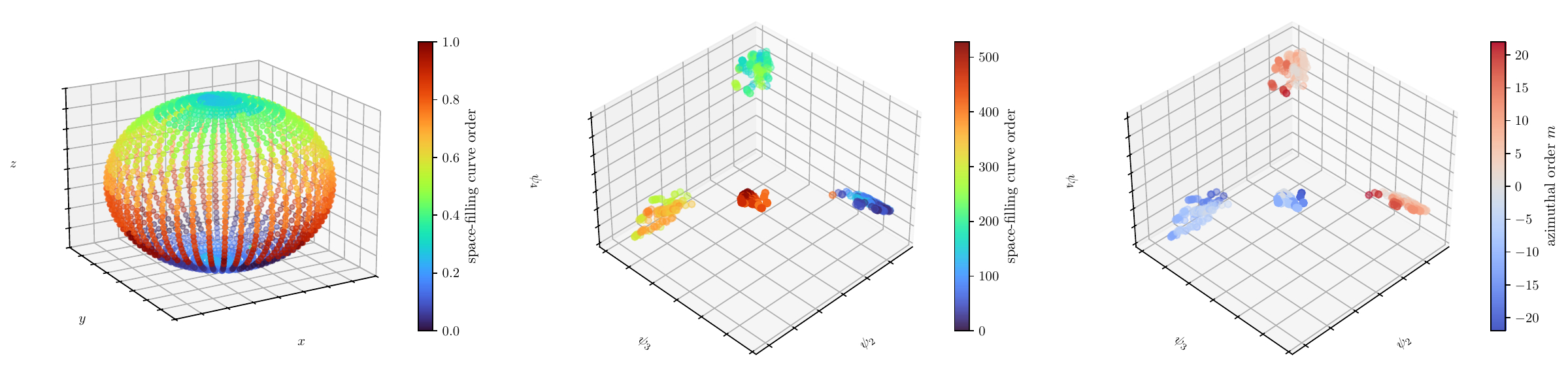}
         \caption{Structured Grid}
         \label{fig:SH_grid_sub}
     \end{subfigure}
     \hfill
     \begin{subfigure}[b]{\linewidth}\centering\includegraphics[width=\linewidth]{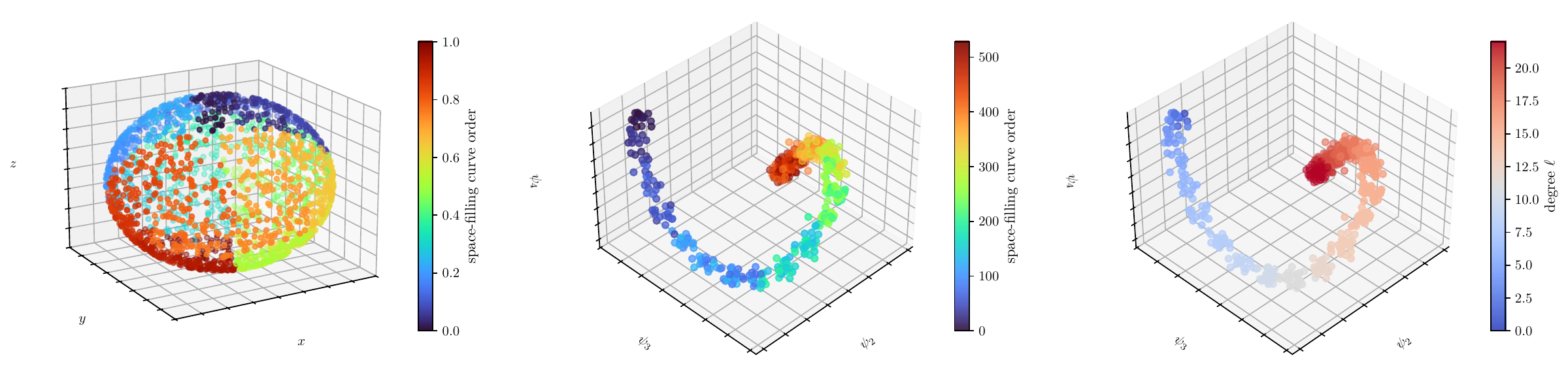}
         \caption{Random Sampling}
         \label{fig:SH_random_sub}
     \end{subfigure}
    \caption{Embedding results for the spherical-harmonic eigenvectors and the sampled points on the unit sphere ($N=2116$, $l_{\max}=23$), colored by the learned space-filling-curve ordering, for (a)~the structured grid and (b)~the random sampling. From left to right in each row: the sampled points on the sphere; the diffusion embedding of the eigenvectors by curve order; and the same embedding colored by frequency, the azimuthal order $m$ for the structured grid and the degree $l$ for the random scheme. The grid resolves the order $m$, whereas the random scheme resolves the degree $l$.}
  
    \label{fig:SH_embedding}
\end{figure}

\begin{figure}
     \centering
     \includegraphics[width=\linewidth]{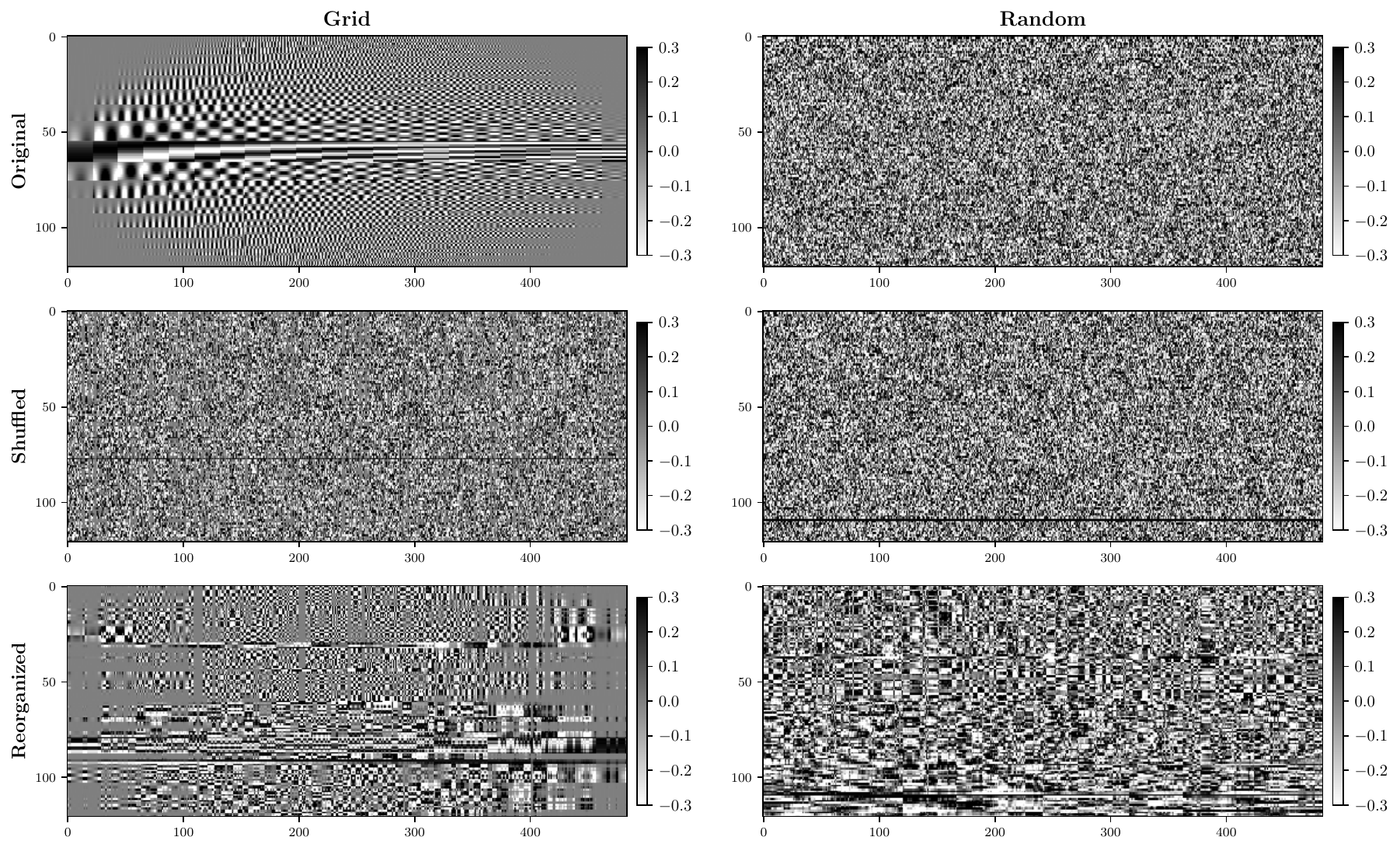}
     \caption{Matrix of spherical-harmonic basis functions for the two sampling schemes (structured grid, left; random sampling, right). Each column shows, from top to bottom, the original matrix, the randomly shuffled matrix, and the matrix reorganized by the flip-perturbed \texttt{Questionnaire}.}
     \label{fig:SH_comparison}
\end{figure}

We apply Algorithm~\ref{alg:questionnaire} to the matrix of spherical-harmonic eigenvectors, with $Y=\{y_j\}$ the eigenfunctions and $X=\{x_i\}$ the sampled points on the sphere. Figure~\ref{fig:SH_embedding} embeds $Y$ (first three eigenvectors of $L_{sym}$ from $W^{corr}_{\mathcal{T}_X}$) and $X$ for a structured grid and a random sampling, colored by the learned curve; Figure~\ref{fig:SH_comparison} shows the $l_{max}=11$ ($N=484$) matrices in original, shuffled, and reorganized states. The grid is already locally regular while the random and permuted matrices look noisy, yet in both schemes the curve recovers the latent block structure.

Because the two schemes sample the sphere differently, the \texttt{Questionnaire} learns different geometries, and hence different space-filling-curve orderings, on both factors. On the structured grid the azimuth is equiangular, so $\cos(m\phi)$ and $\sin(m\phi)$ are exact discrete Fourier modes; the learned embedding of the eigenvectors then clusters the harmonics by azimuthal order $m$, with the degree $l$ varying within each mode. For the random sampling, where the rotation leaves no privileged order $m$, the embedding instead arranges the eigenvectors by degree $l$. Sampled in these two ways, the same operator exposes complementary coordinates of the dual index $(l,m)$.

Table~\ref{table:sh_results} reports the compression results for both schemes across sizes $N$.

\begin{table}[t]
\centering
\makebox[\textwidth][c]{%
\begin{minipage}[t]{0.52\textwidth}
\centering
\subcaption*{Structured grid}
{\scriptsize\renewcommand{\arraystretch}{0.85}\setlength{\tabcolsep}{2pt}%
\resizebox{\linewidth}{!}{%
\begin{tabular}{lcccccc}
    \hline
    \textbf{N} & \textbf{484} & \textbf{1024} & \textbf{2116} & \textbf{4096} & \textbf{8100} & \textbf{16384} \\
    $l_{max}$ & 11 & 14 & 23 & 32& 45&  64 \\ 
    Direct Evaluation &1.30E-4 &5.13E-4 &2.23E-3 &7.64E-3& 3.48E-2 &1.14E-1\\
    \hline
       {\texttt{Questionnaire}}\\
    Prep. Time (sec) &   6.91E1 & 7.77E1 & 9.50E1& 1.54E2 & 2.52E2 & 1.18E3   \\
    \hline
{\texttt{Butterfly}} \\
    Permuted (MiB)  & 9.68E-1 & 5.17E0 & 2.26E1 & 8.47E1 & 2.91E2 & 1.36E3
    \\
    Original (MiB)& 8.02E-1 & 2.46E0 & 6.38E0 & 1.64E1 & 4.97E1 & 9.25E1\\
    Reorganized (MiB)    & 8.59E-1 & 2.56E0 & 6.38E0 & 1.77E1 & 5.80E1 & 9.53E1\\
    
    Prep. Time (sec) & 3.91E-2& 4.04E-1 &3.57E0 &8.19E0& 3.06E1 & 5.77E2    \\
    Forward Transform (sec)   & 1.57E-4 & 3.83E-4& 2.02E-3& 3.86E-3& 1.08E-2& 2.92E-2 \\
    
    Relative $\ell_2$ Error & 5.43E-12& 7.25E-12 & 6.25E-12 & 6.19E-12& 8.58E-12& 1.05E-11  \\ \hline
    {\texttt{eGHWT}}\\
    Permuted (MiB)& 9.70E-1 & 3.93E0& 1.74E1& 6.14E1& 2.40E2  &9.50E2     \\
    Original (MiB)  & 7.01E-1 & 7.83E-1& 1.11E1 & 1.11E1 & 1.27E2 &1.44E2      \\
    Reorganized (MiB) & 4.71E-1 & 9.30E-1 & 6.12E0& 9.41E0 &  6.62E1 & 1.56E2  \\ 
    Prep. Time (sec)    & 2.42E0 & 7.80E0& 3.25E1& 1.48E2& 4.46E2  & 3.64E3\\
     Forward Transform (sec)   & 2.99E-2& 6.60E-2& 2.51E-1& 4.59E-1 & 1.80E0 &  4.16E0    \\
    Relative $\ell_2$ Error    & 4.73E-3& 2.40E-3& 3.92E-3& 2.71E-3 & 3.72E-3& 2.29E-3 \\ \hline
 
    \end{tabular}}}

\vspace{4pt}
\includegraphics[width=\linewidth]{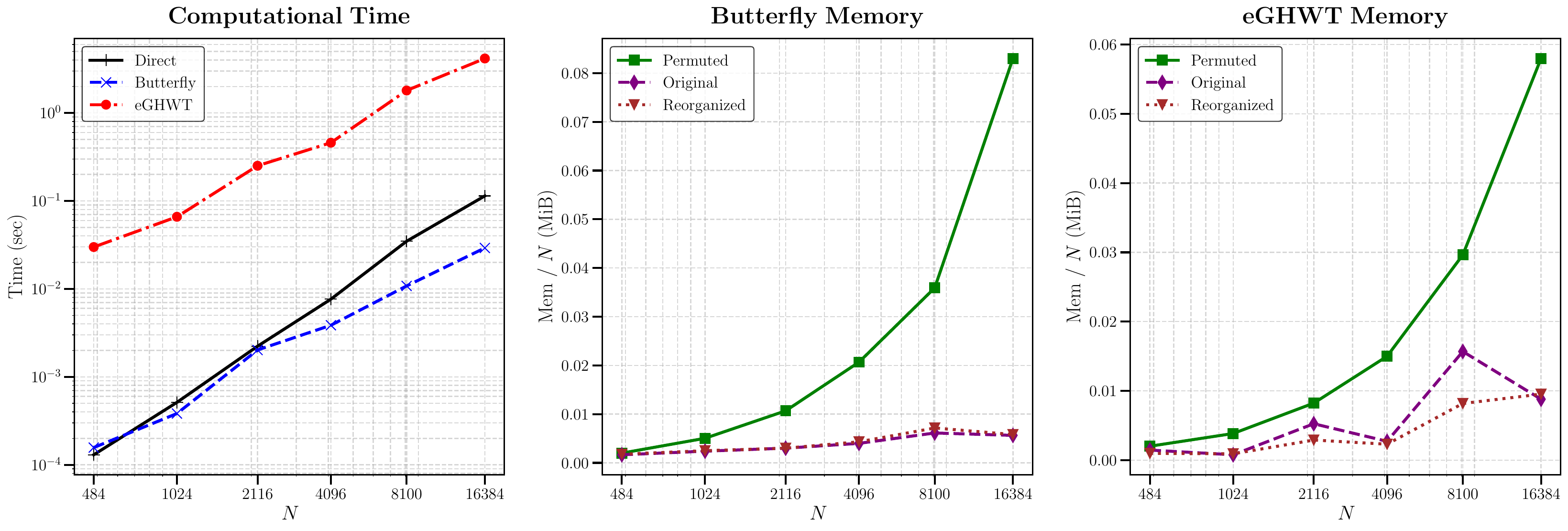}
\end{minipage}\hspace{0.045\textwidth}
\begin{minipage}[t]{0.52\textwidth}
\centering
\subcaption*{Random sampling}
{\scriptsize\renewcommand{\arraystretch}{0.85}\setlength{\tabcolsep}{2pt}%
\resizebox{\linewidth}{!}{%
\begin{tabular}{lcccccc}
    \hline
    \textbf{N} & \textbf{484} & \textbf{1024} & \textbf{2116} & \textbf{4096} & \textbf{8100} & \textbf{16384} \\
    $l_{max}$ & 11 & 14 & 23 & 32 & 45 & 64 \\ 
    Direct Evaluation & 1.05E-4 & 5.08E-4 & 2.23E-3 & 1.13E-2 & 3.91E-2 & 1.34E-1 \\
    \hline
    \multicolumn{7}{l}{\texttt{Questionnaire}} \\
    Prep. Time (sec) & 6.88E1 & 7.72E1 & 9.47E1 & 1.66E2 & 2.69E2 & 1.22E3 \\
    \hline
    \multicolumn{7}{l}{\texttt{Butterfly}} \\
    Permuted (MiB)  & 9.68E-1 & 5.17E0 & 2.26E1 & 8.47E1 & 2.91E2 & 1.36E3 \\
    Original (MiB) & 9.68E-1 & 5.17E0 & 2.26E1 & 8.46E1 & 2.90E2 & 1.33E3 \\ 
    Reorganized (MiB) & 9.68E-1 & 5.06E0 & 1.55E1 & 6.29E1 & 1.60E2 & 4.62E2 \\
    Prep. Time (sec)   & 8.39E-2 & 5.26E-1 & 2.09E0 & 1.01E1 & 2.53E1 & 5.65E2 \\
    Forward Transform (sec)  & 1.93E-4 & 1.42E-3 & 5.06E-3 & 8.09E-3 & 2.04E-2 & 4.13E-2 \\
    Relative $\ell_2$ Error & 5.73E-12 & 6.51E-12 & 9.13E-12 & 1.07E-11 & 3.60E-11 & 1.95E-11 \\ 
    \hline
    \multicolumn{7}{l}{\texttt{eGHWT}} \\
    Permuted (MiB) & 1.07E0 & 4.73E0 & 2.00E1 & 7.39E1 & 2.83E2 & 1.15E3 \\
    Original (MiB) & 1.01E0 & 4.37E0 & 1.82E1 & 6.54E1 & 2.48E2 & 9.71E2 \\
    Reorganized (MiB) & 8.88E-1 & 3.77E0 & 1.55E1 & 5.67E1 & 2.10E2 & 8.30E2 \\ 
    Prep. Time (sec)    & 2.37E0 & 7.80E0 & 3.20E1 & 1.47E2 & 4.46E2 & 3.62E3 \\
    Forward Transform (sec)   & 3.99E-2 & 1.28E-1 & 4.28E-1 & 1.37E0 & 4.84E0 & 1.78E1 \\
    Relative $\ell_2$ Error & 4.73E-3 & 4.50E-3 & 4.80E-3 & 5.00E-3 & 4.72E-3 & 4.94E-3 \\ 
    \hline
    \end{tabular}}}

\vspace{4pt}
\includegraphics[width=\linewidth]{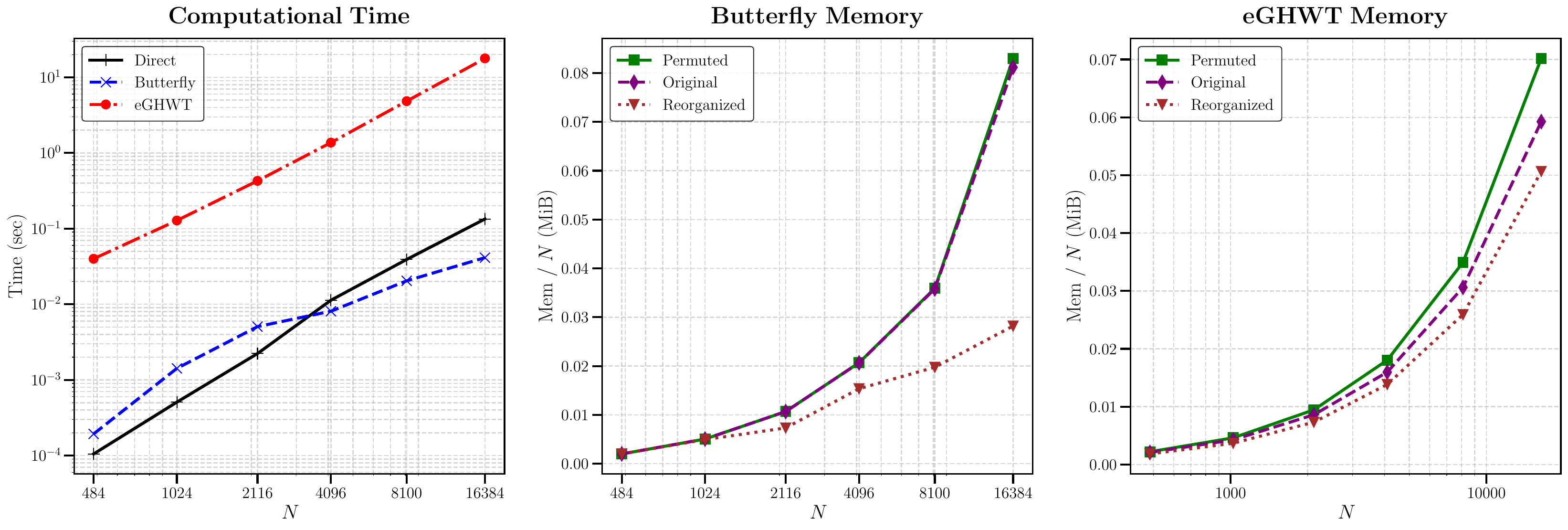}
\end{minipage}}
\caption{Spherical harmonics: times, errors, and memory (tables) and computational time and memory versus $N$ (plots) for the structured grid (left) and random sampling (right).}
\label{table:sh_results}
\end{table}

\subsection{Approximated Laplacian of Riemannian Manifold}\label{sec:manifold}
On a general Riemannian manifold the Laplace--Beltrami eigenfunctions have no closed form and must be approximated numerically; unlike the sphere, such a manifold has no distinguished symmetry group or structured quadrature grid, so its eigenvectors have no known dedicated fast transform and are not a permutation of an operator that does. We build a ``twisted bean'' surface $\mathcal{M}\subset\mathbb{R}^3$ by sampling $N$ points with $(u,v)$ uniform on $[-\tfrac{\pi}{2},\tfrac{\pi}{2}]\times[0,2\pi]$ under the map $\mathbf{x}(u,v)=[x,y,z]^\top$,
$$\begin{aligned}
x &= L \sin u \\
y &= r(u) \cos v \cos(t L \sin u) - \left[ r(u) \sin v + 0.4(L \sin u)^2 \right] \sin(t L \sin u) \\
z &= r(u) \cos v \sin(t L \sin u) + \left[ r(u) \sin v + 0.4(L \sin u)^2 \right] \cos(t L \sin u)
\end{aligned}$$
where $r(u) = W \cos u (1 - 0.2 \cos u)$, with $L = 0.8$, $W = 0.4$, and twist $t = -0.05$, giving a smooth curved manifold of varying local density.

We approximate $\Delta_{\mathcal{M}}$ from a self-tuning Gaussian adjacency
\begin{equation}
W_{ij} = \exp\left( -\frac{\|\mathbf{x}_i - \mathbf{x}_j\|^2}{\sigma_i \sigma_j} \right),
\end{equation}
with $\sigma_i$ the distance to the $30$th-percentile $k$-nearest neighbor, and apply the density-invariant renormalization $W^{(1)} = D^{-1} W D^{-1}$. From the symmetric Laplacian $L_{\mathrm{sym}} = I - D_1^{-1/2} W^{(1)} D_1^{-1/2}$ we solve $L_{\mathrm{sym}} \boldsymbol{\psi} = \lambda \boldsymbol{\psi}$ and map back via $\boldsymbol{\Phi} = D_1^{-1/2} \boldsymbol{\psi}$; the $\lfloor N/4 \rfloor$ smallest-eigenvalue eigenvectors form the basis matrix $\mathcal{K} \in \mathbb{R}^{\lfloor N/4 \rfloor \times N}$.

\begin{figure}[!ht]
\includegraphics[width=\linewidth, trim=1cm 1cm 0cm 1cm, clip]{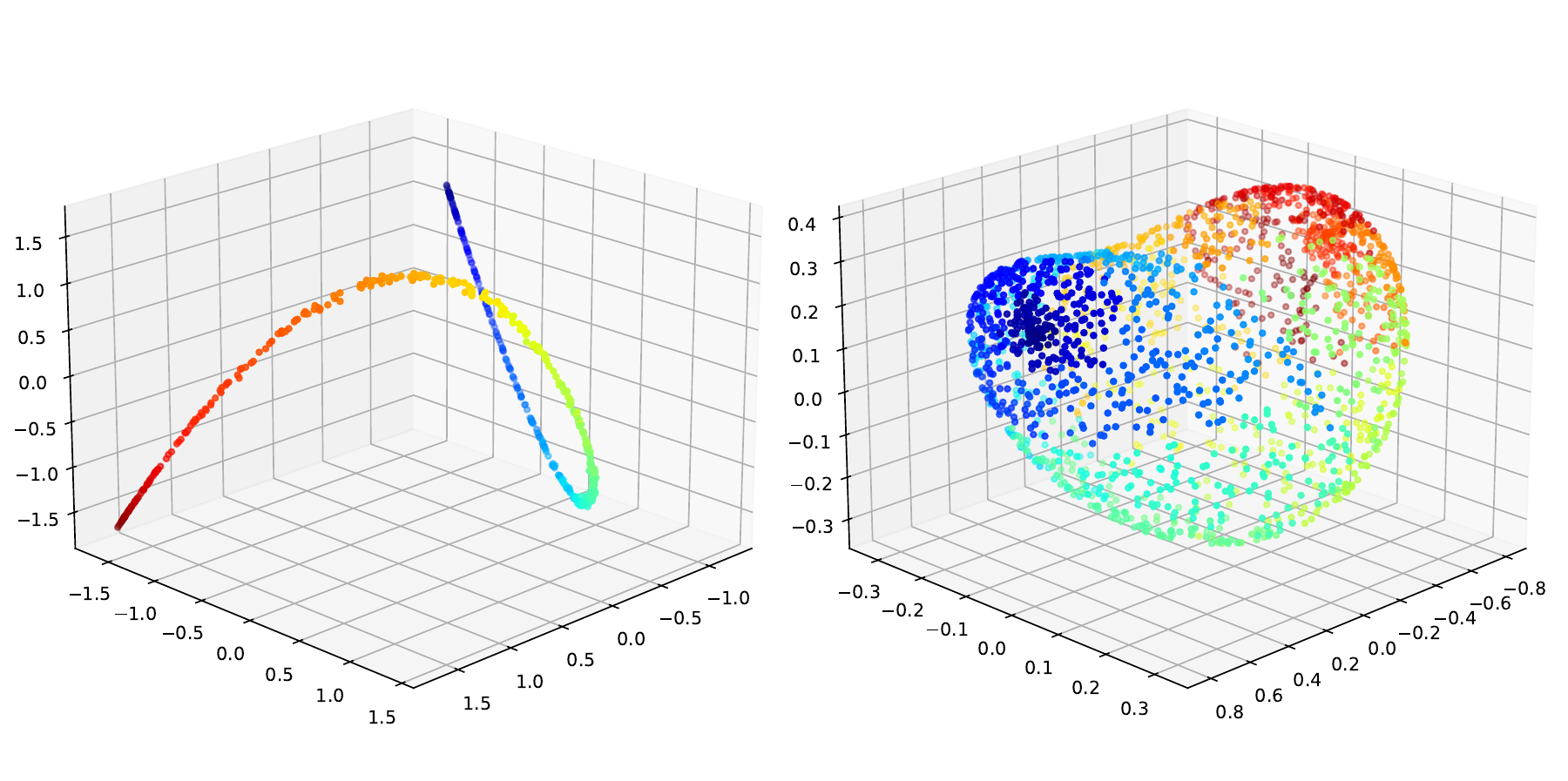}
    \caption{Embedding results for Laplacian eigenvectors and sampled points on the twsited manifold. The color gradient is according to the learned space-filling curve ordering. \textbf{Left:} Diffusion embedding of the eigenvectors.
    \textbf{Right:} Embedding of the sampled points on the twisted manifold.}\label{fig:Bean_embedding}
\end{figure}

\begin{figure}
     \centering
         \includegraphics[width=\linewidth]{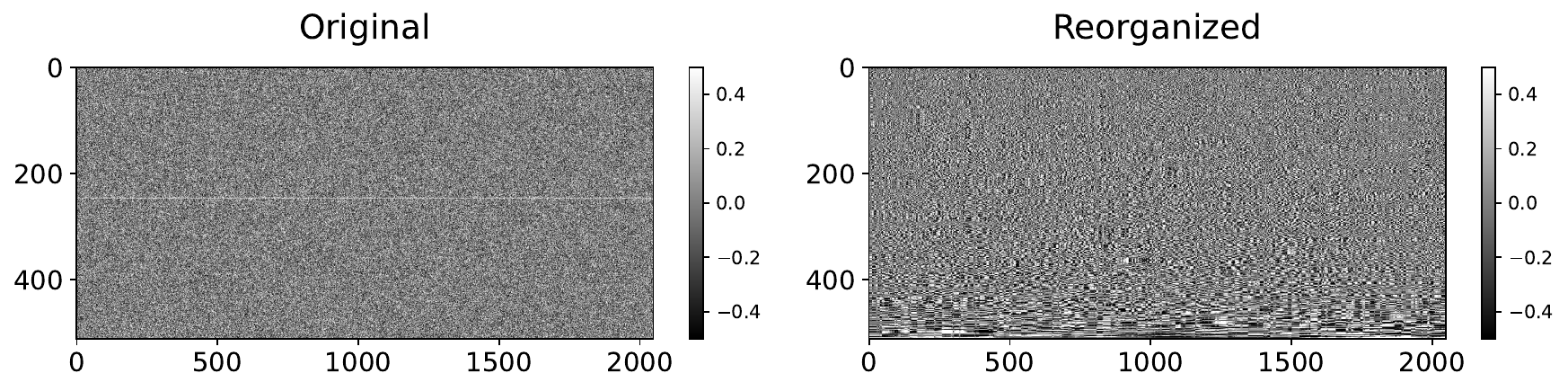}
     \caption{Matrix of Laplacian eigenvectors sampled from the twisted manifold. 
     The left panel displays the randomly permuted matrix, and the right panel displays the reorganized matrix.}
\label{fig:Bean_comparison}
\end{figure}

Figure~\ref{fig:Bean_embedding} embeds the eigenvectors $K$ and the point cloud $X$, colored by the learned curve, Figure~\ref{fig:Bean_comparison} shows the $N=2048$ matrix, noisy under a random permutation and block-structured after reorganization into a picture of the phase space of the manifold, in which points and eigenvectors are in duality, and Table~\ref{table:Bean} reports the compression results across sizes $N$.

\begin{table}[t]
\centering
\makebox[\textwidth][c]{%
\begin{minipage}[c]{0.52\textwidth}
\centering
{\scriptsize\renewcommand{\arraystretch}{0.85}\setlength{\tabcolsep}{2pt}%
\resizebox{\linewidth}{!}{%
\begin{tabular}{lcccccc}
    \hline
    \textbf{N} & \textbf{1024} & \textbf{2048} & \textbf{4096} & \textbf{8192} & \textbf{16384} & \textbf{32768}  \\ 
    Direct Evaluation& 1.85E-4 & 6.54E-4 & 6.81E-3 & 2.61E-2& 9.13E-2 & 3.32E-1\\
    \hline
       {\texttt{Questionnaire}}\\
    Prep. Time (sec) &   8.95E1 & 1.17E2 & 1.82E2& 4.05E2 & 1.75E3 & 5.32E3   \\
    \hline
    {\texttt{Butterfly}} \\
    Permuted (MiB)      & 4.60E0 & 1.73E1 & 8.25E1 & 2.93E2 & 1.35E3 & 4.76E3
    \\
    Reorganized (MiB)   & 3.12E0 & 1.39E1 & 5.30E1 & 1.94E2 & 6.61E2 & 2.14E3   \\
    Prep. Time (sec)      & 5.19E-1  & 1.39E0  & 6.01E0  & 3.96E1  & 1.44E2  & 6.31E2 \\
    Forward Transform (sec)      & 1.56E-4  & 4.52E-4  & 1.56E-3  & 5.20E-3  & 1.61E-2  & 5.17E-2\\
    
    Relative $\ell_2$ Error  & 2.06E-7  & 
    2.18E-7 &  2.07E-7  & 2.05E-7 & 1.92 E-7 & 1.77E-7     \\ \hline
    {\texttt{eGHWT}}\\
    Permuted (MiB)    & 3.43E0     & 1.39E1 & 5.53E1 &2.22E2 & 8.86E2 &-- \\
    Reorganized (MiB)  & 2.68E0& 1.12E1 & 4.48E1& 1.80E2& 7.21E2 & --  \\ 
    Prep. Time (sec)  &7.72E0& 3.68E1& 1.44E2 &  7.23E2&   3.86E3 &--  \\
     Forward Transform (sec) &1.09E-1& 3.46E-1 & 1.18E1  & 4.22E1& 1.57E2 &--    \\
    Relative $\ell_2$ Error    &2.83E-2 & 2.69E-2 & 2.75E-2 & 2.76E-2 & 2.75E-2 &--\\ \hline
    \end{tabular}}}
\end{minipage}\hspace{0.045\textwidth}
\begin{minipage}[c]{0.52\textwidth}
\centering
\includegraphics[width=\linewidth]{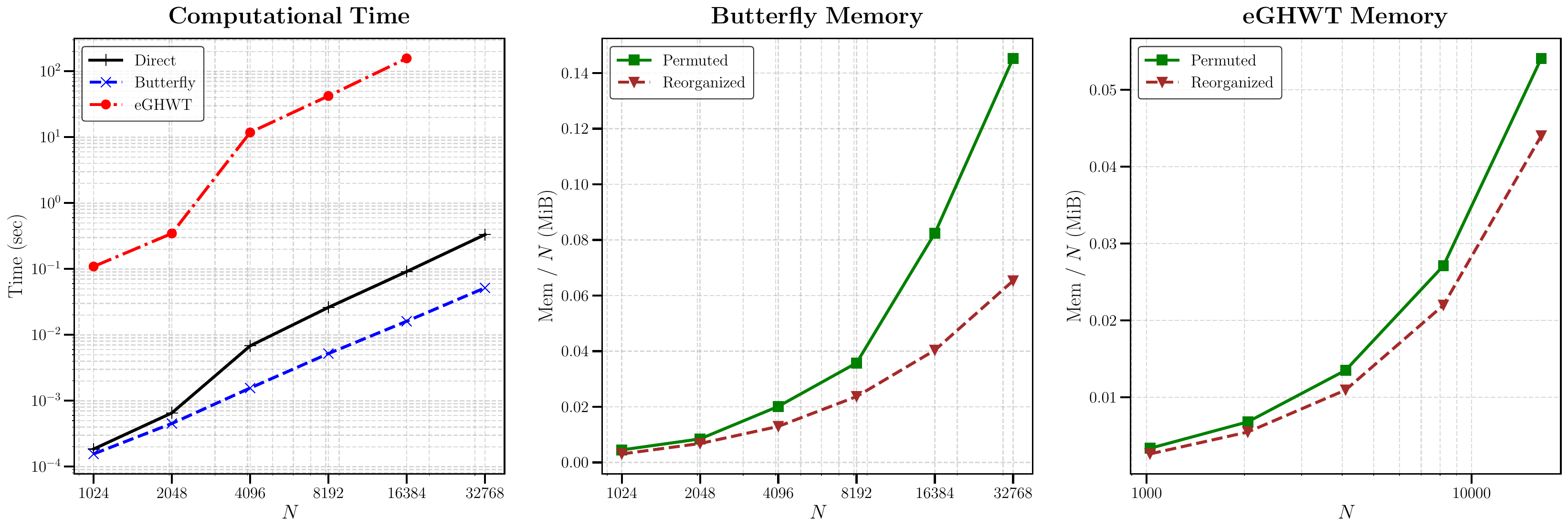}
\end{minipage}}
\caption{Twisted manifold ($L_{rw}$ eigenvectors): times, errors, and memory (left), and computational time and memory versus $N$ (right).}
\label{table:Bean}
\end{table}

\section{Discussion and Future Work}\label{sec:discussion}

We have presented the \texttt{QFFT} framework, which integrates the learning of the dual geometries of the domain and range  of an operator with hierarchical matrix factorization to achieve efficient computation without prior knowledge of the domain's coordinate system. By employing the \texttt{Questionnaire} algorithm, we uncover the latent structures of the row and column spaces, each organized by the geometry learned on the other, so that the reorganized matrix becomes an image  of the operator's phase space: a smooth representation on which the \texttt{Butterfly} and \texttt{eGHWT} algorithms apply.

Our numerical experiments demonstrate that the \texttt{Butterfly} factorization remains highly robust for oscillatory kernels, such as those arising from spherical harmonics and Helmholtz operators, maintaining high precision with $\OO(N \log N)$ complexity. While the \texttt{eGHWT} approach provides significant memory compression for structured and piecewise-smooth operators, its current implementation for the forward transform is limited by single-threaded execution. Future work will focus on the parallelization and algorithmic optimization of the \texttt{eGHWT} transform to make it competitive with direct evaluations at larger scales.

A critical component of this framework is the iterative refinement process within the \texttt{Questionnaire} algorithm. While we have shown empirical success across diverse geometries, the formal theoretical analysis of the convergence criteria and the stability of the fixed points in the dual-tree construction remain open questions. In particular, the convergence of the transform cost is data-dependent, and imposing a smooth weighting in the tree-building step is a natural way to regularize the dual affinities and stabilize the iteration toward a fixed point.
Conceptually we revealed various dual geometries (the fixed points of the questionnaire in a reoganized data matrix). Our aim throughout is to \emph{discover} the geometry under which an $\OO(N\log N)$ computation becomes possible; the \texttt{Butterfly} and \texttt{eGHWT} certify that this geometry is genuinely compressible, while a fast transform calibrated to each discovered geometry, as the FFT is calibrated to the DFT, is a natural next step for future work. We point out that for the example of the matrix of spherical Harmonics on the  3 dimensional sphere we have at least two natural non equivalent phase space organizations corresponding to the  induced Euclidean and  Hermitian structures  in the 4 dimensional ambient space.
In particular, we will explore these mathematical foundations in future work and their extensions to linear and nonlinear higher-dimensional tensor decompositions and more complex transformations  arising in analytic  learning methods.

\section*{Code Availability}
The source code for the proposed method, including the QFFT and examples for DST-IV, Green's functions, Helmholtz kernels, spherical harmonics, and Laplacian kernels, is available at \url{https://github.com/pei-chun-su/QuestionnaireFastTransform}.
\section*{Acknowledgements}
We gratefully acknowledge V. Rokhlin, N. Saito, M. O'Neil, H.W. Zhang, and P. Beckman for their insightful discussions on factorization algorithms and for providing access to their software libraries, which have greatly enhanced both the development of this paper and the associated computational implementations. We also thank R. Lederman,  Tiago Domingues, and the anonymous referees for their constructive discussions, comments, and suggestions.

\appendix
\section{Related work on the alternating construction }

A parallel line to the \texttt{Questionnaire} algorithm casts the row--column coupling as an optimal-transport~\cite{peyre2019computational} fixed-point problem, in which a ground metric on one mode induces a Wasserstein distance on the other and one seeks a pair of metrics that is self-consistent under this map. One method~\cite{huizing2022unsupervised} learns the ground metrics on samples and features simultaneously as the positive singular vectors of the nonlinear map that sends a ground metric to its induced Wasserstein distances, and computes them with a scalable, entropically regularized solver; another~\cite{auffenberg2025unsupervised} recasts the same fixed point as a stochastic random-function iteration; and a further variant~\cite{dusterwald2024fast} accelerates the computation by embedding the two modes on trees and replacing the transport distance with a tree-Wasserstein approximation.

A number of works have sought to establish guarantees on the fixed point of such iterations. The tree-Wasserstein iteration~\cite{lin2026joint} is shown to admit a convergent \emph{subsequence} whose limit is a fixed point when the regularization weights are strictly positive, though convergence of the full iteration from an arbitrary initialization is not guaranteed and the limit is unique only under strong regularization; the optimal-transport formulations establish existence and uniqueness of the fixed point under suitable conditions~\cite{huizing2022unsupervised}, and this iteration~\cite{auffenberg2025unsupervised} was recently shown to converge linearly, relaxing the paracontractivity assumption previously required. These guarantees, however, come with structural restrictions: being built on optimal transport, or on histogram and tree-Wasserstein distances, the variants compare the rows and columns as probability distributions (so the data must be nonnegative and normalized), and they rely on a regularization term~\cite{cuturi2013sinkhorn}, often strong, to force the contractivity-type conditions needed to reach the fixed point. Moreover, these works detect convergence only through the change of the learned affinity or metric from one iteration to the next. Such a criterion certifies that the iteration has stopped moving, but reveals neither \emph{what} geometry it has converged to nor how good that limit is: different initializations may reach different fixed points of differing quality, and a stationarity check can neither tell them apart nor certify that any of them yields a compressible organization, the very property a downstream fast transform requires.

\section{Walsh Wavelet Packets and Best-Basis Selection}\label{app_Walsh}
This section reviews key ideas from the work of Thiele and Villemoes \cite{thiele1996fast}, who introduced a fast, adaptive framework for time–frequency analysis  based on Walsh tilings. 
In view of the fact that we organize the rows or columns of the matrix by an ordered binary tree leading to coordinates of  x or y in their trees, we can parametrize the points by their binaries as points on the interval [0,1]and can use the standard Walsh functions construction as functions on rows or columns as illustrated below for simplicity

\paragraph{Walsh Wavelet Packets} The Walsh system \(\{W_n\}_{n=0}^\infty\) forms a complete orthonormal basis for \(L^2[0,1)\), consisting of piecewise constant functions that take values in \(\{\pm1\}\). These functions are defined recursively by

\begin{align}\label{eq_Walsh}
W_0(t) &= 1, \nonumber\\
W_{2n}(t) &= W_n(2t) + (-1)^n W_n(2t - 1), \nonumber\\
W_{2n+1}(t) &= W_n(2t) - (-1)^n W_n(2t - 1),
\end{align}
for \(t \in [0,1)\), with the convention that \(W_n(t)=0\) for \(t\notin[0,1)\), so that the dilated and translated terms \(W_n(2t)\) and \(W_n(2t-1)\) are well defined; each \(W_n\) has exactly \(n\) sign changes.
A \emph{tile} is a dyadic rectangle of area one, parametrized by scale indices \(j, j', k, n \in \mathbb{N} \cup \{0\}\), with \(k < 2^j\). The associated time interval and sequency band are given by
\begin{equation}\label{eq_tiles}
I := [2^{-j}k, 2^{-j}(k+1)), \quad \omega := [2^{j'} n, 2^{j'} (n+1)).
\end{equation}
When \(j = j'\), the tile \(p = I \times \omega\) is considered unit-area. Let $\mathcal{P}$
be the collection of all of these unit-area dyadic rectangles. A localized Walsh function associated with tile \(p \in \mathcal{P}\) is then defined as
\begin{equation}
w_p(t) = w_{j,k,n}(t) := 2^{j/2} W_n(2^j t - k),
\end{equation}
which is supported on \(I\) and has \(n\) sign changes corresponding to the "sequency" band \(\omega\).

Let \( N = 2^L \) denote the signal length. A discrete signal of length \( N \) can be viewed as a piecewise constant function \( f \colon [0, N) \to \mathbb{R} \), where the function is constant on each unit interval \([n, n+1)\), for \( n = 0, 1, \ldots, N-1 \). We identify this domain \([0,N)\) with the unit interval \([0,1)\) by the rescaling \(t\mapsto t/N\) (equivalently, the space tiles \(I\) above are scaled by \(N\) to \([0,N)\)), so that the atoms \(w_p\) and the signal \(f\) share a common domain and the expansion below is well defined. The Walsh wavelet packet framework constructs an overcomplete dictionary of orthonormal functions by associating each basis atom \( w_p \) with a dyadic space-sequency tile \( p \subset \mathcal{S}_N := [0, N) \times [0, 1) \). At each scale \( j = 0, 1, \ldots, L \), the domain \( \mathcal{S}_N \) is partitioned into \( N \) disjoint tiles, yielding a total of \((L+1)N\) candidate atoms. Any orthonormal basis corresponds to a disjoint partition by tiles \( \mathcal{B} \subset \mathcal{P} \), where \( \mathcal{P} \) denotes the full collection of tiles. The signal admits the expansion
\begin{equation}
f(x) = \sum_{p \in \mathcal{B}} \langle f, w_p \rangle w_p(x),
\end{equation}
where orthogonality holds whenever the tiles are disjoint, i.e., \( p \cap \tilde{p} = \emptyset \) for all \( p \ne \tilde{p} \in \mathcal{B} \). The Walsh transform coefficients \( \langle f, w_p \rangle \) can be computed efficiently in \( \OO(N \log N) \) time using the Fast Walsh–Hadamard Transform (FWHT) \cite{fino1976unified}, which leverages the binary structure of the basis and requires only additions and subtractions. 

This construction extends naturally to two-dimensional signals. Let
\( \mathcal{K} \colon [0, N_Y)\times[0, N_X) \to \mathbb{R} \) denote an image of size \( N_Y \times N_X \), where $N_X = 2^L$ and $N_Y = 2^{L'}$. The two-dimensional Walsh wavelet packet dictionary is formed via tensor products of one-dimensional atoms:
\begin{equation}
w_{(p, q)}(y, x) := w_q(y) w_p(x), \quad (q, p) \in \mathcal{P}_Y \times \mathcal{P}_X,
\end{equation}
where each index \( (q, p) \) corresponds to a dyadic tile in the product domain \( \mathcal{S}_{N_Y} \times \mathcal{S}_{N_X} \). The total number of atoms in this dictionary is \((L'+1)N_Y \times (L+1) N_X\), and any admissible two-dimensional basis corresponds to a tiling \( \mathcal{B} \subset \mathcal{P}_Y \times \mathcal{P}_X \). The image may be represented as
\begin{equation}
\mathcal{K}(y, x) = \sum_{(q, p) \in \mathcal{B}} \langle \mathcal{K}, w_{(q, p)} \rangle w_{(q, p)}(y, x),
\end{equation}
with orthogonality satisfied whenever \( p \cap \tilde{p} = \emptyset \) or \( q \cap \tilde{q} = \emptyset \) for all \( (q, p), (\tilde{q}, \tilde{p}) \in \mathcal{B} \).

\paragraph{Best Basis Selection}
Given a signal \(f \) of length \( N = 2^L \), the objective of the best basis selection is to select a tiling \(\mathcal{B} \subset \mathcal{P}\) such that the associated basis \(\{w_p\}_{p \in \mathcal{B}}\) yields an efficient representation of \(f\). A cost function \(c(p)\) is assigned to each tile \(p \in \mathcal{P}\), typically defined in terms of the magnitude of the transform coefficient \(\langle f, w_p \rangle\). A common choice is
\begin{equation}
c(p) := |\langle f, w_p \rangle|^\ell,
\end{equation}
for \(0 < \ell < 2\). The total cost over a tiling \(\mathcal{B}\) of $\mathcal{S}_N$ is then given by
\begin{equation}
C(f, \mathcal{B}) := \sum_{p \in \mathcal{B}} c(p).
\end{equation}
This cost functional is designed to be small when the energy of \(f\) is concentrated in a small number of significant coefficients \(\langle f, w_p \rangle\), making it well suited for selecting sparse or compact representations.
The \emph{best basis} is the tiling \(\mathcal{B}^* = \argmin_{\mathcal{B}}C(f,\mathcal{B})\), that minimizes the total cost \(C(f, \mathcal{B})\) over all admissible tilings of $\mathcal{S}_N$. Thiele and Villemoes \cite{thiele1996fast} proposed a dynamic programming algorithm that compares the cost of a parent tile with the combined cost of its horizontal and vertical splits, and selects the configuration with the lowest total cost. This recursive process efficiently identifies the optimal decomposition. For a signal with length $N$, the total computational complexity of the method is of $\OO(N \log N)$.

This one-dimensional framework was later extended to two-dimensional signals, particularly in the context of image compression \cite{lindberg2000image}. In the 2D setting, the analysis domain becomes \(\mathcal{K}(y,x)\) on $[0,N_X)\times [0,N_Y)$, and adaptive tilings are constructed over dyadic rectangles in the joint space–sequency plane \(\mathcal{S}_{N_X} \times \mathcal{S}_{N_Y}\). Each 2D tile corresponds to a separable product of 1D tiles. The best-basis selection algorithm applies analogously in two dimensions, recursively minimizing the total cost across hierarchical tilings. The total computational complexity is of $\OO(N_X N_Y(\log N_X\log N_Y))$. We suggest that readers refer to \cite{thiele1996fast, lindberg2000image} for the technical details.

\paragraph{Extension to General Graph Domains}
The Generalized Haar-Walsh Transform (\texttt{GHWT})~\cite{irion2015applied, saito2022eghwt} extends the classical Haar-Walsh wavelet packet framework to signals supported on general graphs. It constructs a hierarchical dictionary of orthonormal basis vectors by recursively applying localized averaging and differencing operations over a tree-structured partition of the vertex set. Unlike classical wavelet packets, which operate on uniformly spaced dyadic intervals with an inherent linear ordering, the \texttt{GHWT} accommodates domains without a global coordinate system. The underlying partition tree is not restricted to binary splits; nodes may be divided into arbitrary subsets, enabling adaptive tilings that reflect the geometry of the underlying graph, while still permitting coefficient computation with \( \OO(N \log N) \) complexity, analogous to the FWHT. In our setting the ``graph'' on which the (e)GHWT operates is simply the row index set \(Y\) (or the column index set \(X\)), equipped with the hierarchical partition tree and dual affinity learned by the \texttt{Questionnaire} (Algorithm~\ref{alg:questionnaire_basic}); the (e)GHWT best basis is then built on that learned tree, which is how the transform connects to the problem of this paper.

The extended \texttt{GHWT} (\texttt{eGHWT})~\cite{saito2022eghwt} further enhances this framework by implementing a best-basis selection algorithm over the \texttt{GHWT} dictionary, analogous to the one- and two-dimensional best-basis algorithms developed for Haar-Walsh wavelet packets. This allows for the adaptive selection of an orthonormal basis that minimizes a user-defined cost functional, such as entropy or sparsity.
The algorithm is applicable to signals of arbitrary length on graphs of general structure, without requiring \( N = 2^L \). The overall computational complexity is \( \OO(N \log N) \) for one-dimensional signals or graph-supported data under unbalanced hierarchical partitions. For two-dimensional signals, the complexity becomes $\OO(N_X N_Y(\log N_X\log N_Y))$, analogous to the classical 2D setting. We refer the reader to~\cite{irion2015applied, saito2022eghwt} for algorithmic details.

\paragraph{Ordering Ambiguity}
The construction of a hierarchical tree is central to the \texttt{eGHWT}, as it defines how basis vectors are recursively partitioned and labeled. However, in many cases, especially when no underlying coordinate system is available for the data domains \( X \) and \( Y \), the structure of the tree is not uniquely determined. This ambiguity arises not only in the choice of partition boundaries, but also in the labeling of subregions (e.g., “left” vs. “right”) at each level of the tree, which affects the binary tags assigned to basis functions. As a result, the outcome of \texttt{eGHWT} best-basis selection can vary, leading to different basis representations and transform coefficients for the same underlying signal and data.
In~\cite{saito2022eghwt}, the tree is constructed using eigenvectors of the graph Laplacian to recursively partition the row and column domains. We introduce a finer and more stable approach to tree construction in such coordinate-free settings (see Algorithm~\ref{alg:affinity_ordering}).

\paragraph{Fast \texttt{eGHWT} Kernel Application}
Once the reorganized kernel \( \mathcal{K}(x, y) \) has been compressed (Algorithm~\ref{alg:compression}), it is represented in the best-basis form
\[
w_{(q_i, p_i)}(x, y) = w_{q_i}(x) w_{p_i}(y),
\]
indexed by \( (q_i, p_i) \in \mathcal{B}^* \), along with the corresponding \texttt{GHWT} coefficients \( c_{(q_i, p_i)} \) for \( i = 1, \ldots, N_Y N_X \). To reduce storage and computational cost, only the top \( n \) coefficients are retained, which sufficiently satisfy a prescribed error \( \varepsilon \).
Given an input vector \( f \in \mathbb{R}^{N_X} \), interpreted as a function on the discrete domain \( [0, N_X) \), we compute its \texttt{GHWT} expansion coefficients:
\[
c_p^f := \langle f, w_p \rangle,
\]
which can be obtained in \( \OO(N_X \log N_X) \) time. The matrix-vector multiplication is then performed in the transform domain by aggregating over the retained basis pairs:
\[
\widehat{c}_{q_i} := \sum_{\substack{p_i \text{ such that } \\ (q_i, p_i) \in \mathcal{B}^*}} c_{(q_i, p_i)} \cdot c_{p_i}^f, \quad i = 1, \ldots, n,
\]
followed by reconstruction of the output vector via the inverse \texttt{eGHWT}:
\begin{equation}\label{eq_GHWT_apply}
\widehat{\mathcal{K}}f(y) := \sum_{\substack{p_i \text{ such that } \\ (q_i, p_i) \in \mathcal{B}^*}} \widehat{c}_{q_i} \cdot w_{q_i}(y).
\end{equation}
When \( n \) is of \(\OO(N \log N) \) or smaller, the total complexity of this procedure is \( \OO(N \log N) \), offering a significant computational advantage over the \( \OO(N^2) \) cost of direct matrix-vector multiplication.

\section{Interpolative Decomposition}\label{app_ID}
\begin{lemma} (Interpolative Decomposition~\cite{cheng2005compression,martinsson2007interpolation}).
Let $A \in \mathbb{R}^{m \times n}$ and $k$ be a positive integer such that $k \le \min(m, n)$. There exists a matrix $B \in \mathbb{R}^{m \times k}$ consisting of a subset of the columns of $A$, and a matrix $P \in \mathbb{R}^{k \times n}$, such that:
\begin{enumerate}
    \item $P$ contains a $k \times k$ identity matrix as a submatrix.
    \item $\max_{i,j} |P_{ij}| \le 2$.
    \item The spectral norm satisfies $\|P\|_2 \le \sqrt{4k(n-k) + 1}$.
    \item The $k$-th singular value satisfies $\sigma_k(P) \ge 1$.
    \item If $k = \min(m, n)$, then $A = BP$.
    \item If $k < \min(m, n)$, the approximation error is bounded by:
    \[
    \|A - BP\|_2 \le \sqrt{4k(n-k) + 1} \, \sigma_{k+1}(A)
    \]
\end{enumerate}
\end{lemma}

\section{Butterfly Factorization}\label{app_Butterfly}

The butterfly algorithm is a multiscale fast transform designed to efficiently apply dense integral operators arising in special function transforms. These operators yield kernel matrices that are globally dense but exhibit a localized low-rank structure when restricted to sub-blocks corresponding to well-separated regions in the input and output domains.  We are guided by the prototypical oscillatory Fourier kernel
$\mathcal{K}(y, x) = e^{i \gamma x y / 4}$,
 The result of~\cite{landau1962prolate, o2010algorithm} shows that the numerical rank of the associated integral operator, to a precision \( \varepsilon \), is controlled by both the frequency parameter \( \gamma \) and the size of the domain $\times$ range to which the kernel is restricted.
Similar results for the Fourier--Bessel integral operator are also  derived in the same paper~\cite{o2010algorithm}. 

The butterfly algorithm leverages the localized low-rank structure of many oscillatory kernels. Given a kernel matrix \( \mathcal{K}(y, x) \) with \( x \in X \), \( y \in Y \), and \( |X| = |Y| = N \), the algorithm begins by constructing dyadic partition trees \( \mathcal{T}_X \) and \( \mathcal{T}_Y \) in the output and input domains, respectively. These trees define a multilevel block structure, where each level \( \ell \in \{0, 1, \dots, L\} \), with \( L = \log_2 N \), corresponds to interactions between node pairs at the level \( \ell \) of \( \mathcal{T}_X \) and the level \( L - \ell \) of \( \mathcal{T}_Y \). In practice the dyadic subdivision is stopped before reaching singletons: it descends only until every leaf folder \( V_k^\ell \) holds at most \( \max_k |V_k^\ell| \) columns, and skeletonization begins from these leaf blocks. This maximum leaf-block size is a user-chosen parameter of the method that fixes the finest scale of the factorization, trading preparation cost against the tightness of the compression.
At each level \( \ell \), the algorithm applies interpolative decomposition (\texttt{ID})~\cite{cheng2005compression,goreinov2001maximal,martinsson2007interpolation} to approximate the sub-matrices corresponding to these block interactions. Specifically, given a matrix \( A \in \mathbb{C}^{m \times n} \) whose \((r+1)\)st singular value is sufficiently small, \texttt{ID} yields an approximation
\begin{equation}
A \approx BP.
\end{equation}

\begin{figure}
    \centering
\includegraphics[width=1.05\linewidth]{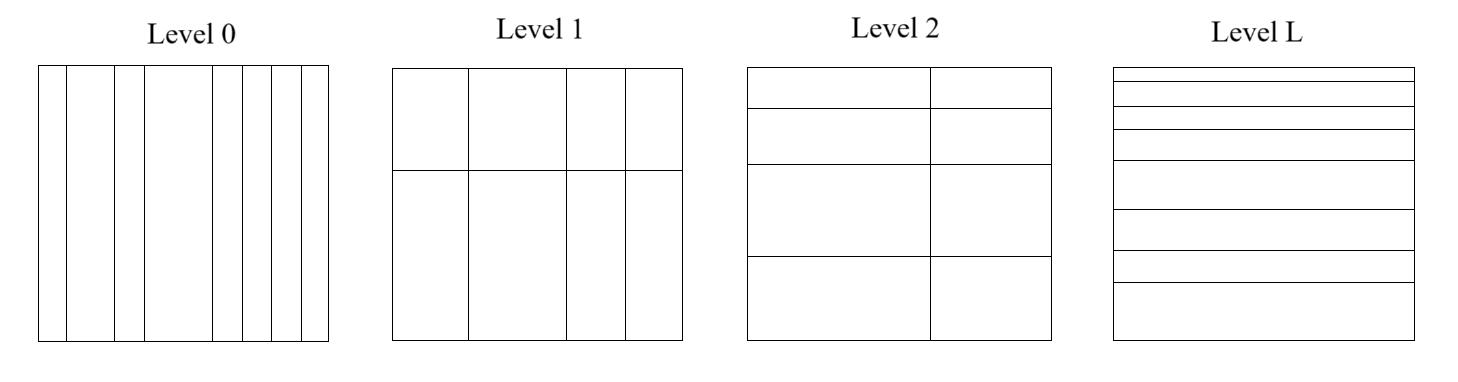}
    \caption{Skeletonization}
    \label{fig:skeletonization}
\end{figure}

Let two adjacent blocks at level \( \ell \) be approximated by interpolative decompositions as
\begin{equation}
B_l^{\ell} = B_{l}^{\ell+1} P^{\ell}_l, \quad B_r^{\ell} = B_{r}^{\ell+1} P^{\ell}_r.
\end{equation}
At level \( \ell +1\), $B^{\ell+1}_l$ and $B^{\ell+1}_r$ merge to form new blocks. As illustrated in Figure~\ref{fig:skeletonization}, the upper portion of the parent block is constructed by merging the upper halves of \( B_l^{\ell+1} \) and \( B_r^{\ell+1} \), while the lower part merges the lower halves. The process of selecting representative rows or columns during this merge is referred to as \emph{skeletonization}.
This recursive divide-and-merge strategy continues until the coarsest level \( L \) is reached, resulting in a multiscale hierarchical factorization of the original matrix. Assuming that each \texttt{ID} is performed with error $\varepsilon$  ,the final factorization enables efficient application of the matrix to a vector up to precision $\varepsilon$:
\begin{equation}\label{eq_bf_apply}
(\mathcal{K}f)(y_i) = 
\sum_{j_1} B^{(L)}_{i, j_{L-1}}
\sum_{j_2} P^{(L-1)}_{j_{L-1}, j_{L-2}}
\cdots
\sum_{j} P^{(0)}_{j_{0}, j} f(x_{j}) 
+ \OO(\varepsilon),
\end{equation}
where each interpolation matrix \( P^{(\ell)} \) maps the function values at skeleton points at level \( \ell \) to approximations on larger blocks at level \( \ell + 1 \), and \( B^{(L)} \) reconstructs the final output values using the skeletons at the finest scale. 

The precomputation of the \texttt{Butterfly} factorization requires \( \OO\big((N\log N)^2\big) \) operations due to the need to perform \texttt{ID} on all bi-folder submatrices of the kernel. Assuming that the numerical rank \( r \) of the merged blocks at each level remains comparable to that of the blocks at the previous level, the resulting factorization stores interpolation matrices and the basis block of the coarsest-level \( B^{(L)} \), with total memory complexity \( \OO(rN \log N) \). As for the kernel $\mathcal{K}(y, x) = e^{i \gamma x y / 4}$, this rank-stability assumption typically holds for kernels when the row and column domains are similarly structured at each stage of the hierarchy.
In \cite{o2010algorithm, tygert2010fast}, the method is applied to a broad class of oscillatory and special-function kernels. Once constructed, the matrix-vector multiplication can be performed in \( \OO(N \log N) \) time. For a detailed discussion of the precomputation procedure, hierarchical compression, and fast matrix-vector application, we refer the reader to~\cite{o2010algorithm, tygert2010fast}.

\bibliographystyle{plain}
\bibliography{ev}    

\end{document}